\documentclass[11pt]{amsart}
\usepackage{amsmath,amssymb,wasysym}
\newtheorem{theorem}{Theorem}[section]
\newtheorem{corollary}[theorem]{Corollary}
\newtheorem{lemma}[theorem]{Lemma}
\newtheorem{proposition}[theorem]{Proposition}

\theoremstyle{definition}
\newtheorem{definition}[theorem]{Definition}
\newtheorem{remark}[theorem]{Remark}

\begin{document}

\title[A fully nonlinear flow for two-convex hypersurfaces]{A fully nonlinear flow for two-convex hypersurfaces in Riemannian manifolds}
\author{Simon Brendle and Gerhard Huisken}
\address{Department of Mathematics \\ Columbia University \\ 2990 Broadway \\ New York, NY 10027}
\address{Mathematisches Institut \\ Universit\"at T\"ubingen \\ 72076 T\"ubingen \\ Germany}
\begin{abstract}
We consider a one-parameter family of closed, embedded hypersurfaces moving with normal velocity $G_\kappa = \big ( \sum_{i < j} \frac{1}{\lambda_i+\lambda_j-2\kappa} \big )^{-1}$, where $\lambda_1 \leq \hdots \leq \lambda_n$ denote the curvature eigenvalues and $\kappa$ is a nonnegative constant. This defines a fully nonlinear parabolic equation, provided that $\lambda_1+\lambda_2>2\kappa$. In contrast to mean curvature flow, this flow preserves the condition $\lambda_1+\lambda_2>2\kappa$ in a general ambient manifold. 

Our main goal in this paper is to extend the surgery algorithm of Huisken-Sinestrari to this fully nonlinear flow. This is the first construction of this kind for a fully nonlinear flow. As a corollary, we show that a compact Riemannian manifold satisfying $\overline{R}_{1313}+\overline{R}_{2323} \geq -2\kappa^2$ with non-empty boundary satisfying $\lambda_1+\lambda_2 > 2\kappa$ is diffeomorphic to a $1$-handlebody.

The main technical advance is the pointwise curvature derivative estimate. The proof of this estimate requires a new argument, as the existing techniques for mean curvature flow due to Huisken-Sinestrari, Haslhofer-Kleiner, and Brian White cannot be generalized to the fully nonlinear setting. To establish this estimate, we employ an induction-on-scales argument; this relies on a combination of several ingredients, including the almost convexity estimate, the inscribed radius estimate, as well as a regularity result for radial graphs. We expect that this technique will be useful in other situations as well.
\end{abstract} 
\maketitle

\section{Introduction}

Throughout this paper, we fix an integer $n \geq 3$ and a real number $\kappa \geq 0$. We consider a closed, embedded hypersurface $M_0$ in an $(n+1)$-dimensional Riemannian manifold which is $\kappa$-two-convex in the sense that $\lambda_1+\lambda_2 > 2\kappa$. We evolve $M_0$ with normal velocity 
\[G_\kappa = \Big ( \sum_{i < j} \frac{1}{\lambda_i+\lambda_j-2\kappa} \Big )^{-1},\] 
where $\lambda_1 \leq \hdots \leq \lambda_n$ denote the principal curvatures. This defines a fully nonlinear parabolic evolution equation. The case $\kappa=0$ is particularly interesting. In this case, we require that the hypersurface $M_0$ is two-convex, and we evolve $M_0$ with normal velocity 
\[G = \Big ( \sum_{i < j} \frac{1}{\lambda_i+\lambda_j} \Big )^{-1}.\] 
In the first part of this paper, we analyze the properties of this flow up to the first singular time: 

\begin{theorem}
\label{analysis.up.to.first.singular.time}
Let $M_t = \partial \Omega_t$, $t \in [0,T)$, be a one-parameter family of closed, embedded, $\kappa$-two-convex hypersurfaces in a compact Riemannian manifold which move with velocity $G_\kappa$. Then the following statements hold: \\
(i) The function $G_\kappa$ is uniformly bounded from below on bounded time intervals. Moreover, if the curvature tensor of the ambient manifold satisfies $\overline{R}_{1313}+\overline{R}_{2323} \geq -2\kappa^2$ at each point on $M_t$, then $\inf_{M_t} G_\kappa$ blows up in finite time. \\
(ii) The ratio $\frac{\lambda_1+\lambda_2-2\kappa}{H}$ is uniformly bounded from below on bounded time intervals. \\
(iii) The hypersurfaces $M_t$ are almost convex at points where the curvature is large. More precisely, given $\delta>0$, we can find a positive constant $K$, depending only on $\delta$, $T$, $\kappa$, the initial hypersurface $M_0$, and the ambient manifold, such that $\lambda_1 \geq -\delta G_\kappa$ whenever $G_\kappa \geq K$. \\
(iv) Given $\delta>0$, we can find positive constants $\eta$ and $K$, depending only on $\delta$, $T$, $\kappa$, the initial hypersurface $M_0$, and the ambient manifold, such that $\lambda_n-\lambda_2 \leq \delta G_\kappa$ whenever $G_\kappa \geq K$ and $\lambda_1 \leq \eta G_\kappa$. \\
(v) At each point on $M_t$, the inscribed radius is bounded from below by $\frac{\alpha}{G_\kappa}$. Here, $\alpha$ is a positive constant that depends only on $T$, $\kappa$, the initial hypersurface $M_0$, and the ambient manifold. \\ 
(vi) The quantity $G_\kappa^{-2} \, |\nabla h| + G_\kappa^{-3} \, |\nabla^2 h|$ is uniformly bounded from above at all points where the curvature is sufficienly large. Again, the constants depend only on $T$, $\kappa$, the initial hypersurface $M_0$, and the ambient manifold.
\end{theorem} 

We note that the curvature condition $\overline{R}_{1313}+\overline{R}_{2323} \geq -2\kappa^2$ in statement (i) is sharp. Indeed, if $\kappa<1$ and the ambient manifold is a standard hyperbolic cusp, then there exists a family of hypersurfaces moving with speed $G_\kappa$ which exists for all $t \in [0,\infty)$.

The statement (i) follows easily from the maximum principle. Moreover, (ii) is a straightforward adaptation of results due to Andrews \cite{Andrews1}. The proof will be described in Section \ref{basics}. 

The statements (iii) and (iv) are consequences of Theorem \ref{cylindrical.estimate} below. The proof of the cylindrical estimate uses the Michael-Simon Sobolev inequality and Stampacchia iteration. This is discussed in Section \ref{cyl.est}. 

The statement (v) was established by Andrews-Langford-McCoy \cite{Andrews-Langford-McCoy1} when $\kappa=0$ and the ambient manifold is the Euclidean space $\mathbb{R}^{n+1}$. Their work easily carries over to the case $\kappa \geq 0$. In the Riemannian setting, various error terms arise due to the background geometry, but these can be controlled in the same way as in \cite{Brendle4}. We note that the corresponding noncollapsing estimate for embedded, mean convex solutions of mean curvature flow was first established in the fundamental work of Brian White \cite{White1},\cite{White2}. In \cite{Sheng-Wang}, Sheng and Wang gave an alternative proof of the noncollapsing estimate for mean curvature flow. Later, Andrews \cite{Andrews2} gave another proof of that estimate based on a direct maximum principle argument. Finally, in \cite{Brendle3} and \cite{Brendle4}, the first author improved this to a sharp estimate: more precisely, for an embedded, mean convex solution of mean curvature flow, the inscribed radius is bounded from below by $\frac{1-\delta}{H}$ at points where the curvature is large. We note that a similar estimate holds for the fully nonlinear flow considered in this paper (cf. \cite{Brendle-Hung}), but we will not use this stronger estimate here. The arguments in \cite{Andrews2},\cite{Andrews-Langford-McCoy1}, and \cite{Brendle3} are based on an application of the maximum principle to a suitably chosen function that depends on a pair of points. This technique originated in earlier work of the second author \cite{Huisken} on the curve shortening flow in the plane (see also \cite{Grayson},\cite{Hamilton2}). A recent survey can be found in \cite{Brendle2}. 

The pointwise curvature derivative estimate in statement (vi) is the most difficult part of Theorem \ref{analysis.up.to.first.singular.time}. The corresponding estimate for mean curvature flow was established by Brian White \cite{White1},\cite{White2} in the mean convex case (see also \cite{Haslhofer-Kleiner1} and \cite{Sheng-Wang}), and by the second author and Carlo Sinestrari \cite{Huisken-Sinestrari2} under the stronger assumption of two-convexity. The arguments in \cite{Haslhofer-Kleiner1},\cite{White1},\cite{White2} rely on the monotonicity formula for mean curvature flow, whereas the proof in Section 6 of \cite{Huisken-Sinestrari2} is based on the maximum principle. The fully nonlinear case requires a new argument, as there is no analogue of the monotonicity formula and a direct maximum principle argument does not seem to work. In the following, we sketch the main ideas that allow us to overcome this obstacle. Let us consider a point $(\bar{x},\bar{t})$ in spacetime where the curvature is very large. Using the inscribed radius estimate, we can find a point $p$ such that the ball $B_{\alpha G_\kappa(\bar{x},\bar{t})^{-1}}(p)$ lies inside $M_{\bar{t}}$ and touches $M_{\bar{t}}$ at $\bar{x}$. Given any point $x \in B_{2\alpha G_\kappa(\bar{x},\bar{t})^{-1}}(p)$, we construct a pseudo-cone $C_{p,x}$, which has a conical singularity at $x$ with some fixed opening angle. In geodesic normal coordinates around $x$, the boundary of $C_{p,x}$ is a rotationally symmetric hypersurface with the property that the curvature in radial direction is bounded from above by a small negative multiple of $d(p,x)^{-1}$. We then distinguish two cases:

Suppose first that the hypersurfaces $M_t$ can be represented as radial graphs in a parabolic neighborhood of the point $(\bar{x},\bar{t})$ with size comparable to $G_\kappa(\bar{x},\bar{t})^{-1}$. In this case, a regularity result for radial graphs (cf. Section \ref{rad.graph}) gives an upper bound for $G_\kappa^{-2} \, |\nabla h| + G_\kappa^{-3} \, |\nabla^2 h|$ at the point $(\bar{x},\bar{t})$.

Suppose next that the hypersurfaces $M_t$ cannot be represented as radial graphs in a suitable parabolic neighborhood of the point $(\bar{x},\bar{t})$. In this case, we can find a time $\tilde{t} \leq \bar{t}$ and a point $\tilde{x} \in B_{2\alpha G_\kappa(\bar{x},\bar{t})^{-1}}(p)$ with the property that the pseudo-cone $C_{p,\tilde{x}}$ lies inside $M_{\tilde{t}}$ and touches $M_{\tilde{t}}$ from the inside at some point $y \in M_{\tilde{t}}$. Since the radial curvature of the pseudo-cone is bounded from above by a negative multiple of $G_\kappa(\bar{x},\bar{t})$, it follows that $\frac{\lambda_1(y,\tilde{t})}{G_\kappa(\bar{x},\bar{t})}$ is bounded from above by a negative constant. The almost convexity property in statement (iii) then implies that $G_\kappa(y,\tilde{t})$ is much larger than $G_\kappa(\bar{x},\bar{t})$. We now invoke the Neck Detection Lemma to conclude that the point $y$ lies at the center of a neck which is contained in $M_{\tilde{t}}$. Since the pseudo-cone $C_{p,\tilde{x}}$ lies inside $M_{\tilde{t}}$, this setup contradicts elementary geometry.

To summarize, we are able to prove the curvature derivative estimate, assuming that the Neck Detection Lemma can be applied. However, the proof of the Neck Detection Lemma relies in a crucial way on the pointwise curvature derivative estimate! To avoid a circular argument, we observe that, in order to prove the curvature derivative estimate at $(\bar{x},\bar{t})$, we need to apply the Neck Detection Lemma at $(y,\tilde{t})$, and the curvature at $(y,\tilde{t})$ is much larger than the curvature at $(\bar{x},\bar{t})$. This allows us to carry out an induction-on-scales argument. The details are discussed in Section \ref{derivative.est}. 

In the second part of this paper, we use a surgery procedure as in \cite{Huisken-Sinestrari2} to extend the flow beyond singularities. 

\begin{theorem}
\label{existence.of.surgically.modified.flows}
Let $M_0 = \partial \Omega_0$ be a closed, embedded, $\kappa$-two-convex hypersurface in a compact Riemannian manifold. Given any $T>0$, there exists a surgically modified flow with velocity $G_\kappa$ which starts from $M_0$ and is defined on the time interval $[0,T)$. Moreover, if the curvature tensor of the ambient manifold satisfies $\overline{R}_{1313}+\overline{R}_{2323} \geq -2\kappa^2$ at each point in $\Omega_0$, then the flow becomes extinct in finite time.
\end{theorem}

As a consequence of Theorem \ref{existence.of.surgically.modified.flows}, we obtain the following classification of diffeomorphism types (see also \cite{Sha}):

\begin{corollary}
A compact Riemannian manifold satisfying $\overline{R}_{1313}+\overline{R}_{2323} \geq -2\kappa^2$ with non-empty boundary satisfying $\lambda_1+\lambda_2 > 2\kappa$ is diffeomorphic to a $1$-handlebody.
\end{corollary}

The idea of extending solutions of geometric flows past singularities by means of a surgery procedure goes back to the groundbreaking work of Richard Hamilton \cite{Hamilton3},\cite{Hamilton4} on the formation of singularities in the Ricci flow. In particular, in \cite{Hamilton4}, Hamilton developed a surgery algorithm for the Ricci flow on four-manifolds with positive isotropic curvature. In a spectacular series of papers \cite{Perelman1},\cite{Perelman2},\cite{Perelman3}, Perelman successfully implemented a surgery algorithm for the Ricci flow in dimension $3$, and used it to prove the Poincar\'e and Geometrization Conjectures. In \cite{Huisken-Sinestrari2}, the second author and Carlo Sinestrari introduced a notion of mean curvature flow with surgery for two-convex hypersurfaces in Euclidean space $\mathbb{R}^{n+1}$, where $n \geq 3$. The remaining case $n=2$ was recently settled by the authors in \cite{Brendle-Huisken1},\cite{Brendle-Huisken2}; an alternative construction was given by Haslhofer and Kleiner \cite{Haslhofer-Kleiner2}. Unlike Theorem \ref{existence.of.surgically.modified.flows}, the main result in \cite{Huisken-Sinestrari2} cannot be extended to hypersurfaces in a Riemannian manifold: indeed, a two-convex hypersurface in Riemannian manifold may not remain two-convex when evolved by the mean curvature flow.

The proof of Theorem \ref{existence.of.surgically.modified.flows} is presented in Section \ref{estimates.for.flows.with.surgery} and Section \ref{surgery.riemannian.case}. In Section \ref{estimates.for.flows.with.surgery}, we show that the a-priori estimates in Theorem \ref{analysis.up.to.first.singular.time} still hold for surgically modified flows. These a-priori estimates enable us to implement the surgery algorithm from \cite{Huisken-Sinestrari2}. This is completely straightforward if the ambient manifold is the Euclidean space $\mathbb{R}^{n+1}$. Indeed, having established the convexity estimate, the cylindrical estimate, and the curvature derivative estimate for surgically modified flows, the arguments in Section 7 and Section 8 of \cite{Huisken-Sinestrari2} (in particular, the Neck Detection Lemma, the Neck Continuation Theorem, and the surgery algorithm) carry over unchanged to our situation. Finally, extending the results in Section 7 and Section 8 of \cite{Huisken-Sinestrari2} to the Riemannian setting requires some minor adaptations; these are explained in Section \ref{surgery.riemannian.case} below.

\begin{remark}
The exact choice of the normal velocity $G_\kappa$ is not very important. All we need is that $G_\kappa$ satisfies the following structure conditions: 
\begin{itemize}
\item $G_\kappa$ is smooth positive function which is defined on the set of all symmetric matrices satisfying $\lambda_1+\lambda_2 > 2\kappa$. Moreover, $G_\kappa$ approaches $0$ on the boundary of that set. 
\item $G_\kappa$ is a homogeneous function of degree $1$ in $\lambda_1-\kappa,\hdots,\lambda_n-\kappa$. 
\item We have $0 \leq \frac{d}{ds} G_\kappa(h+sA) \big |_{s=0} \leq C \, \text{\rm tr}(A)$ whenever $A$ is two-nonnegative. Moreover, the inequalities are strict unless $A=0$.
\item We have $\frac{d^2}{ds^2} G_\kappa(h+sA) \big |_{s=0} \leq 0$ for every symmetric matrix $A$. Moreover, the inequality is strict unless $A$ is a scalar multiple of $h-\kappa g$.
\end{itemize}
\end{remark}

\textbf{Acknowledgments.} We would like to thank Connor Mooney and Xu-Jia Wang for discussions. We are very grateful to Richard Hamilton for discussions on the non-conic estimate for the Ricci flow. The first author is grateful to Columbia University, the Fields Institute, Toronto, and T\"ubingen University, where parts of this work were carried out. This project was supported by the National Science Foundation under grants DMS-1201924 and DMS-1505724.

\section{Basic properties}

\label{basics}

In this section, we establish some basic properties of the fully nonlinear flow defined above. First, we observe that $G_\kappa$ depends smoothly on the components of $h$; this is a consequence of Theorem 5.7 in \cite{Ball}. Moreover, we clearly have $G_\kappa \leq C(n) \, (H-n\kappa)$ and $\frac{\partial G_\kappa}{\partial h_{ij}} \leq C(n) \, g_{ij}$, where $C(n)$ is a positive constant that depends only on the dimension. We next compute the second derivatives of $G_\kappa$ with respect to $h$.

\begin{proposition}
\label{2nd.variation.of.G}
Suppose that $h$ and $A$ are symmetric $n \times n$ matrices, and that $h$ satisfies $\lambda_1+\lambda_2 > 2\kappa$. Then 
\begin{align*} 
&\frac{d^2}{ds^2} G_\kappa(h+sA) \Big |_{s=0} \\ 
&= -G_\kappa^2 \sum_{i \neq l} \sum_{j \notin \{i,l\}} \frac{1}{(\lambda_i+\lambda_j-2\kappa)(\lambda_l+\lambda_j-2\kappa)} \, \Big ( \frac{1}{\lambda_i+\lambda_j-2\kappa}+\frac{1}{\lambda_l+\lambda_j-2\kappa} \Big ) \, A_{il}^2 \\ 
&- 2 \, G_\kappa^2 \sum_{i<j} \frac{1}{(\lambda_i+\lambda_j-2\kappa)^3} \, (A_{ii}+A_{jj})^2 \\ 
&+ 2 \, G_\kappa^3 \, \Big ( \sum_{i<j} \frac{1}{(\lambda_i+\lambda_j-2\kappa)^2} \, (A_{ii}+A_{jj}) \Big )^2, 
\end{align*} 
where $\lambda_1 \leq \hdots \leq \lambda_n$ denote the eigenvalues of $h$ and $e_1,\hdots,e_n$ are the corresponding eigenvectors.
\end{proposition}

\textbf{Proof.} 
Straightforward calculation. \\

\begin{corollary} 
\label{concavity}
Suppose that $h$ and $A$ are symmetric $n \times n$ matrices, and that $h$ satisfies $\lambda_1+\lambda_2 > 2\kappa$. Then $\frac{d^2}{ds^2} G_\kappa(h+sA) \big |_{s=0} \leq 0$, and equality holds if and only if $A$ is a scalar multiple of $h-\kappa g$. 
\end{corollary}

\textbf{Proof.} 
The inequality $\frac{d^2}{ds^2} G_\kappa(h+sA) \big |_{s=0} \leq 0$ follows immediately from Proposition \ref{2nd.variation.of.G}. Suppose next that equality holds. Then $A_{ij} = 0$ for $i \neq j$. Moreover, we have $A_{ii}+A_{jj} = a \, (\lambda_i+\lambda_j-2\kappa)$ for $i \neq j$, where $a$ is a real number which does not depend on $i$ and $j$. This implies that $A$ is a scalar multiple of $h-\kappa g$. \\

Let $M_t$ be a one-parameter family of closed, embedded, $\kappa$-two-convex hypersurfaces in an $(n+1)$-dimensional compact Riemannian manifold $X$. We assume that the hypersurfaces $M_t$ move inward with normal velocity 
\[G_\kappa = \Big ( \sum_{i < j} \frac{1}{\lambda_i+\lambda_j-2\kappa} \Big )^{-1},\] 
where $\lambda_1,\hdots,\lambda_n$ denote the principal curvatures. The evolution equation of $G_\kappa$ is 
\[\frac{\partial}{\partial t} G_\kappa = \frac{\partial G_\kappa}{\partial h_i^j} \, \frac{\partial}{\partial t} h_i^j = \frac{\partial G_\kappa}{\partial h_{ij}} \, (D_i D_j G_\kappa + h_{ik} \, h_{jk} \, G_\kappa + \overline{R}_{i\nu j\nu} \, G_\kappa).\] 
In the remainder of this section, we discuss two basic a-priori estimates. First, we establish a lower bound for $G_\kappa$; this estimate is needed to ensure that the flow becomes extinct in finite time. Second, we prove that, on any given bounded time interval, the mean curvature is bounded from above by a constant multiple of $G_\kappa$. Both estimates are easy adaptations of Theorem 4.1 in \cite{Andrews1}. 

\begin{lemma}
\label{lower.bound.for.G}
We have $G_\kappa \geq \frac{1}{C} \, e^{-Ct}$, where $C$ is a positive constant that depends only on $\kappa$, the initial hypersurface $M_0$, and the ambient manifold. Moreover, if the curvature tensor of the ambient manifold satisfies $\overline{R}_{1313}+\overline{R}_{2323} \geq -2\kappa^2$ at each point on $M_t$, then $\inf_{M_t} G_\kappa$ approaches infinity in finite time.
\end{lemma} 

\textbf{Proof.} 
Recall that $\frac{\partial G_\kappa}{\partial h_{ij}} \leq C(n) \, g_{ij}$. This implies 
\[\frac{\partial}{\partial t} G_\kappa \geq \frac{\partial G_\kappa}{\partial h_{ij}} \, D_i D_j G_\kappa - C \, G_\kappa.\] 
Using the maximum principle, we obtain $G_\kappa \geq \frac{1}{C} \, e^{-Ct}$, where $C$ is a large constant which is independent of $t$. 

We now assume that the curvature tensor of the ambient manifold satisfies $\overline{R}_{1313}+\overline{R}_{2323} \geq -2\kappa^2$ at each point on $M_t$. Using the identity $G_\kappa = \frac{\partial G_\kappa}{\partial h_{ij}} \, (h_{ij} - \kappa g_{ij})$, we obtain 
\begin{align*} 
\frac{\partial}{\partial t} G_\kappa 
&= \frac{\partial G_\kappa}{\partial h_{ij}} \, D_i D_j G_\kappa + \frac{\partial G_\kappa}{\partial h_{ij}} \, (h_{ik} - \kappa g_{ik}) \, (h_{jk}-\kappa g_{jk}) \, G_\kappa \\ 
&+ \frac{\partial G_\kappa}{\partial h_{ij}} \, (\overline{R}_{i\nu j\nu}+\kappa^2 g_{ij}) + 2\kappa \, G_\kappa^2. 
\end{align*}
Moreover, it follows from the Cauchy-Schwarz inequality that 
\begin{align*} 
G_\kappa^2 
&= \Big ( \sum_{i,j} \frac{\partial G_\kappa}{\partial h_{ij}} \, (h_{ij}-\kappa g_{ij}) \Big )^2 \\ 
&\leq \Big ( \sum_{i,j} \frac{\partial G_\kappa}{\partial h_{ij}} \, g_{ij} \Big ) \, \Big ( \sum_{i,j,k} \frac{\partial G_\kappa}{\partial h_{ij}} \, (h_{ik}-\kappa g_{ik}) \, (h_{jk}-\kappa g_{jk}) \Big ) \\ 
&\leq C(n) \, \sum_{i,j,k} \frac{\partial G_\kappa}{\partial h_{ij}} \, (h_{ik}-\kappa g_{ik}) \, (h_{jk}-\kappa g_{jk}). 
\end{align*} 
Finally, our assumption on the sectional curvature of the ambient manifold implies that the tensor $\overline{R}_{i\nu j\nu}+\kappa^2 g_{ij}$ is two-nonnegative. This implies $\frac{\partial G_\kappa}{\partial h_{ij}} \, (\overline{R}_{i\nu j\nu}+\kappa^2 g_{ij}) \geq 0$. Putting these facts together, we obtain 
\[\frac{\partial}{\partial t} G_\kappa \geq \frac{\partial G_\kappa}{\partial h_{ij}} \, D_i D_j G_\kappa + \frac{1}{C(n)} \, G_\kappa^3.\] 
Using the maximum principle, we conclude that $\inf_{M_t} G_\kappa$ approaches infinity in finite time. This completes the proof of Lemma \ref{lower.bound.for.G}. \\

We next recall the evolution equation for the mean curvature from \cite{Andrews1}. Using the inequality $G_\kappa \leq \frac{\partial G_\kappa}{\partial h_{ij}} \, h_{ij}$, we obtain 
\begin{align*} 
\frac{\partial}{\partial t} H 
&= \Delta G_\kappa + |h|^2 \, G_\kappa + \overline{\text{\rm Ric}}(\nu,\nu)\, G_\kappa \\ 
&\leq \frac{\partial G_\kappa}{\partial h_{ij}} \, (\Delta h_{ij} + |h|^2 \, h_{ij}) + \frac{\partial^2 G_\kappa}{\partial h_{ij} \, \partial h_{kl}} \, D_p h_{ij} \, D_p h_{kl} + \overline{\text{\rm Ric}}(\nu,\nu) \, G_\kappa \\ 
&\leq \frac{\partial G_\kappa}{\partial h_{ij}} \, (D_i D_j H + h_{ik} \, h_{jk} \, H) + \frac{\partial^2 G_\kappa}{\partial h_{ij} \, \partial h_{kl}} \, D_p h_{ij} \, D_p h_{kl} + C \, H + C, 
\end{align*} 
where $C$ is a positive constant that depends only on the ambient manifold. As in \cite{Andrews1}, this evolution equation implies that $\frac{H}{G_\kappa}$ is bounded from above:

\begin{proposition}
\label{uniformly.two.convex}
We have $G_\kappa \geq \beta H$ for all $t \in [0,T)$, where $\beta$ is a positive constant that depends only on $T$, $\kappa$, the initial hypersurface $M_0$, and the ambient manifold. In particular, the ratio $\frac{\lambda_1+\lambda_2-2\kappa}{H}$ is uniformly bounded from below on any bounded time interval. 
\end{proposition}

\textbf{Proof.} 
Recall that 
\[\frac{\partial}{\partial t} G_\kappa \geq \frac{\partial G_\kappa}{\partial h_{ij}} \, (D_i D_j G_\kappa + h_{ik} \, h_{jk} \, G_\kappa) - C \, G_\kappa.\] 
By Corollary \ref{concavity}, $G_\kappa$ is a concave function of the second fundamental form. This implies 
\[\frac{\partial}{\partial t} H \leq \frac{\partial G_\kappa}{\partial h_{ij}} \, (D_i D_j H + h_{ik} \, h_{jk} \, H) + C \, H + C.\] 
Moreover, by Lemma \ref{lower.bound.for.G}, we have $G_\kappa \geq \frac{1}{C}$ for some positive constant $C$ that depends only on $T$, $\kappa$, the initial hypersurface $M_0$, and the ambient manifold. This implies $H \geq \frac{1}{C}$, hence 
\[\frac{\partial}{\partial t} H \leq \frac{\partial G_\kappa}{\partial h_{ij}} \, (D_i D_j H + h_{ik} \, h_{jk} \, H) + C \, H.\] 
Using the maximum principle, we conclude that $\frac{H}{G_\kappa} \leq C$, where $C$ is a positive constant that depends only on $T$, $\kappa$, the initial hypersurface $M_0$, and the ambient manifold. This completes the proof of Proposition \ref{uniformly.two.convex}. \\

Proposition \ref{uniformly.two.convex} implies that 
\[\frac{\partial G_\kappa}{\partial h_{ij}} \geq \frac{G_\kappa^2}{H^2} \, g_{ij} \geq \beta^2 \, g_{ij}.\] 
Therefore, the equation is uniformly parabolic.

\section{The cylindrical estimate}

\label{cyl.est}

Our goal in this section is to prove the following cylindrical estimate:

\begin{theorem}[Cylindrical Estimate]  
\label{cylindrical.estimate}
Let $M_t$, $t \in [0,T)$, be a family of closed, $\kappa$-two-convex hypersurfaces moving with speed $G_\kappa$, and let $\delta$ be an arbitrary positive real number. Then 
\[H \leq \frac{(n-1)^2(n+2)}{4} \, (1+\delta) \, G_\kappa + C,\] 
where $C$ is a positive constant that depends only on $\delta$, $T$, $\kappa$, the initial hypersurface $M_0$, and the ambient manifold.
\end{theorem}
 
In the following, we describe the proof of Theorem \ref{cylindrical.estimate}. In the following lemma, we combine the evolution equation for the mean curvature with the strict concavity property established in Corollary \ref{concavity}.

\begin{lemma} 
\label{good.term}
We have
\[\frac{\partial}{\partial t} H \leq \frac{\partial G_\kappa}{\partial h_{ij}} \, (D_i D_j H + h_{ik} \, h_{jk} \, H) - \frac{1}{C} \, \frac{|\nabla h|^2}{G_\kappa} + C \, G_\kappa\] 
for all $t \in [0,T)$. Here, $C$ is a positive constant that depends only on $T$, $\kappa$, the initial hypersurface $M_0$, and the ambient manifold.
\end{lemma}

\textbf{Proof.} 
Recall that $G_\kappa \geq \beta H$ by Proposition \ref{uniformly.two.convex}. Using Corollary \ref{concavity}, we obtain 
\[\sum_{i,j,k,l} \frac{\partial^2 G_\kappa}{\partial h_{ij} \, \partial h_{kl}} \, A_{ij} \, A_{kl} \leq -\frac{1}{C \, (H-n\kappa)} \, \Big | A - \frac{\text{\rm tr}(A)}{H-n\kappa} \, (h-\kappa g) \Big |^2,\] 
where $C$ is a positive constant that depends on the constant $\beta$ from Proposition \ref{uniformly.two.convex}. This implies 
\[\sum_{i,j,k,l,p} \frac{\partial^2 G_\kappa}{\partial h_{ij} \, \partial h_{kl}} \, D_p h_{ij} \, D_p h_{kl} \leq -\frac{1}{C \, (H-n\kappa)} \, \sum_{i,j,p} \Big ( D_p h_{ij} - \frac{D_p H}{H-n\kappa} \, (h_{ij}-\kappa g_{ij}) \Big )^2.\] 
Using the Codazzi equations, we obtain 
\begin{align*} 
|\nabla H|^2 
&\leq C \sum_{i,j,p} \Big ( -\frac{D_p H}{H-n\kappa} \, (h_{ij}-\kappa g_{ij}) + \frac{D_i H}{H-n\kappa} \, (h_{pj}-\kappa g_{pj}) \Big )^2 \\ 
&\leq C \sum_{i,j,p} \Big ( D_p h_{ij} - \frac{D_p H}{H-n\kappa} \, (h_{ij}-\kappa g_{ij}) - D_i h_{pj} + \frac{D_i H}{H-n\kappa} \, (h_{pj}-\kappa g_{pj}) \Big )^2 + C \\ 
&\leq C \sum_{i,j,p} \Big ( D_p h_{ij} - \frac{D_p H}{H-n\kappa} \, (h_{ij}-\kappa g_{ij}) \Big )^2 + C, 
\end{align*} 
hence 
\[|\nabla h|^2 \leq C \sum_{i,j,p} \Big ( D_p h_{ij} - \frac{D_p H}{H-n\kappa} \, (h_{ij}-\kappa g_{ij}) \Big )^2 + C.\] 
Putting these facts together, we conclude that 
\[\sum_{i,j,k,l,p} \frac{\partial^2 G_\kappa}{\partial h_{ij} \, \partial h_{kl}} \, D_p h_{ij} \, D_p h_{kl} \leq -\frac{1}{C} \, \frac{|\nabla h|^2}{H-n\kappa} + \frac{C}{H-n\kappa},\] 
where $C$ is a positive constant that depends only on $T$, $\kappa$, the initial hypersurface $M_0$, and the ambient manifold. Substituting this into the evolution equation for $H$ gives 
\[\frac{\partial}{\partial t} H \leq \frac{\partial G_\kappa}{\partial h_{ij}} \, (D_i D_j H + h_{ik} \, h_{jk} \, H) - \frac{1}{C} \, \frac{|\nabla h|^2}{H-n\kappa} + C \, H + C + \frac{C}{H-n\kappa}.\] From this, the assertion follows easily. \\

In the following, we fix a positive number $\delta>0$. For $\sigma \in (0,\frac{1}{2})$, we define
\[f_\sigma = G_\kappa^{\sigma-1} \, \Big ( H - \frac{(n-1)^2(n+2)}{4} \, (1+\delta) \, G_\kappa \Big )\] 
and 
\[f_{\sigma,+} = \max \{f_\sigma,0\}.\] 

\begin{proposition}
\label{Lp.bound}
Given any $\delta>0$, we can find a positive constant $c_0$, depending only on $\delta$, $T$, $\kappa$, the initial hypersurface $M_0$, and the ambient manifold, with the following property: if $p \geq \frac{1}{c_0}$ and $\sigma \leq c_0 \, p^{-\frac{1}{2}}$, then we have 
\[\frac{d}{dt} \bigg ( \int_{M_t} f_{\sigma,+}^p \bigg ) \leq (Cp)^p \, |M_t|.\] 
Here, $C$ is a positive constant that depends only on $\delta$, $T$, $\kappa$, the initial hypersurface $M_0$, and the ambient manifold, but not on $\sigma$ and $p$.
\end{proposition}

\textbf{Proof.} 
Using Lemma \ref{good.term}, we obtain 
\begin{align*} 
&\frac{\partial}{\partial t} f_\sigma - \frac{\partial G_\kappa}{\partial h_{ij}} \, D_i D_j f_\sigma - 2(1-\sigma) \, \frac{\partial G_\kappa}{\partial h_{ij}} \, \frac{D_i G_\kappa}{G_\kappa} \, D_j f_\sigma - \sigma \, f_\sigma \, \frac{\partial G_\kappa}{\partial h_{ij}} \, h_{ik} \, h_{jk} \\ 
&+ \sigma (1-\sigma) \, f_\sigma \, \frac{\partial G_\kappa}{\partial h_{ij}} \, \frac{D_i G_\kappa}{G_\kappa} \, \frac{D_j G_\kappa}{G_\kappa} \\ 
&= \Big ( (\sigma-1) \, G_\kappa^{\sigma-2} \, H - \frac{(n-1)^2(n+2)}{4} \, (1+\delta) \, \sigma \, G_\kappa^{\sigma-1} \Big ) \, \Big ( \frac{\partial}{\partial t} G_\kappa - \frac{\partial G_\kappa}{\partial h_{ij}} \, D_i D_j G_\kappa \Big ) \\ 
&+ G_\kappa^{\sigma-1} \, \Big ( \frac{\partial}{\partial t} H - \frac{\partial G_\kappa}{\partial h_{ij}} \, D_i D_j H \Big ) - \sigma \, f_\sigma \, \frac{\partial G_\kappa}{\partial h_{ij}} \, h_{ik} \, h_{jk} \\ 
&\leq -\frac{1}{C} \, G_\kappa^{\sigma-2} \, |\nabla h|^2 + C \, G_\kappa^\sigma. 
\end{align*} 
Consequently, we have 
\begin{align*} 
\frac{\partial}{\partial t} f_\sigma 
&\leq \frac{\partial G_\kappa}{\partial h_{ij}} \, D_i D_j f_\sigma + C \, \frac{|\nabla G_\kappa|}{G_\kappa} \, |\nabla f_\sigma| \\ 
&+ C \, \sigma \, G_\kappa^2 \, f_\sigma - \frac{1}{C} \, G_\kappa^{\sigma-2} \, |\nabla h|^2 + C \, G_\kappa^\sigma 
\end{align*} 
on the set $\{f_\sigma \geq 0\}$. This implies  
\begin{align*} 
\frac{d}{dt} \bigg ( \int_{M_t} f_{\sigma,+}^p \bigg ) 
&\leq p \int_{M_t} f_{\sigma,+}^{p-1} \, \frac{\partial G_\kappa}{\partial h_{ij}} \, D_i D_j f_\sigma + Cp \int_{M_t} f_{\sigma,+}^{p-1} \, \frac{|\nabla G_\kappa|}{G_\kappa} \, |\nabla f_\sigma| \\ 
&+ C \sigma p \int_{M_t} G_\kappa^2 \, f_{\sigma,+}^p - \frac{1}{C} \, p \int_{M_t} G_\kappa^{\sigma-2} \, f_{\sigma,+}^{p-1} \, |\nabla h|^2 \\ 
&+ Cp \int_{M_t} G_\kappa^\sigma \, f_{\sigma,+}^{p-1} - \frac{1}{C} \int_{M_t} G_\kappa^2 \, f_{\sigma,+}^p, 
\end{align*} 
where the last term arises due to the change of the measure. Integration by parts gives 
\begin{align*} 
\frac{d}{dt} \bigg ( \int_{M_t} f_{\sigma,+}^p \bigg ) 
&\leq -p(p-1) \int_{M_t} f_{\sigma,+}^{p-2} \, \frac{\partial G_\kappa}{\partial h_{ij}} \, D_i f_\sigma \, D_j f_\sigma + Cp \int_{M_t} f_{\sigma,+}^{p-1} \, \frac{|\nabla h|}{G_\kappa} \, |\nabla f_\sigma| \\ 
&+ C \sigma p \int_{M_t} G_\kappa^2 \, f_{\sigma,+}^p - \frac{1}{C} \, p \int_{M_t} f_{\sigma,+}^p \, \frac{|\nabla h|^2}{G_\kappa^2} \\ 
&+ Cp \int_{M_t} G_\kappa^\sigma \, f_{\sigma,+}^{p-1} - \frac{1}{C} \int_{M_t} G_\kappa^2 \, f_{\sigma,+}^p, 
\end{align*} 
hence 
\begin{align*} 
\frac{d}{dt} \bigg ( \int_{M_t} f_{\sigma,+}^p \bigg ) 
&\leq -\frac{1}{C} \, p(p-1) \int_{M_t} f_{\sigma,+}^{p-2} \, |\nabla f_\sigma|^2 - \frac{1}{C} \, p \int_{M_t} f_{\sigma,+}^p \, \frac{|\nabla h|^2}{G_\kappa^2} \\ 
&+ C \sigma p \int_{M_t} G_\kappa^2 \, f_{\sigma,+}^p + Cp \int_{M_t} G_\kappa^\sigma \, f_{\sigma,+}^{p-1} - \frac{1}{C} \int_{M_t} G_\kappa^2 \, f_{\sigma,+}^p 
\end{align*} 
for $p$ sufficiently large. To estimate the term $\int_{M_t} G_\kappa^2 \, f_{\sigma,+}^p$, we consider the tensor 
\[S_{ijkl} = -h_{ik} \, h_{jp} \, h_{pl} + h_{jk} \, h_{ip} \, h_{pl} - h_{il} \, h_{jp} \, h_{pk} + h_{jl} \, h_{ip} \, h_{pk}.\] 
A standard commutator identity gives 
\[|D_i D_j h_{kl} - D_j D_i h_{kl} + S_{ijkl}| \leq C \, |h|,\] 
where $C$ depends only on the ambient manifold. This implies 
\begin{align*} 
\int_{M_t} \frac{f_{\sigma,+}^p}{G_\kappa^4} \, |S|^2 
&\leq -\int_{M_t} \frac{f_{\sigma,+}^p}{G_\kappa^4} \, S_{ijkl} \, (D_i D_j h_{kl} - D_j D_i h_{kl}) + C \int_{M_t} \frac{f_{\sigma,+}^p}{G_\kappa^4} \, |S| \, |h| \\ 
&\leq Cp \int_{M_t} f_{\sigma,+}^{p-1} \, \frac{|\nabla h|}{G_\kappa} \, |\nabla f_\sigma| + C \int_{M_t} f_{\sigma,+}^p \, \frac{|\nabla h|^2}{G_\kappa^2} + C \int_{M_t} f_{\sigma,+}^p. 
\end{align*} 
In the next step, we will estimate $|S|^2$ from below. If we diagonalize $h$, then we obtain $S_{ijij} = \lambda_i\lambda_j (\lambda_i-\lambda_j)$ for $i \neq j$. Thus, $|S|^2 \geq \sum_{i,j=1}^n \lambda_i^2\lambda_j^2 (\lambda_i-\lambda_j)^2$. Hence, if $H \geq \frac{(n-1)^2(n+2)}{4} \, (1+\delta) \, G_\kappa$, then $|S|^2 \geq \sum_{i,j=1}^n \lambda_i^2\lambda_j^2 (\lambda_i-\lambda_j)^2 \geq \frac{1}{C} \, G_\kappa^6 - C$ for some constant $C$ which depends only on $\delta$, $T$, $\kappa$, the initial hypersurface $M_0$, and the ambient manifold. In particular, we have $|S|^2 \geq \frac{1}{C} \, G_\kappa^6 - C$ on the set $\{f_\sigma \geq 0\}$. Thus, we conclude that 
\[\int_{M_t} G_\kappa^2 \, f_{\sigma,+}^p \leq Cp \int_{M_t} f_{\sigma,+}^{p-1} \, \frac{|\nabla h|}{G_\kappa} \, |\nabla f_\sigma| + C \int_{M_t} f_{\sigma,+}^p \, \frac{|\nabla h|^2}{G_\kappa^2} + C \int_{M_t} f_{\sigma,+}^p,\] 
where $C$ depends only on $\delta$, $T$, $\kappa$, the initial hypersurface $M_0$, and the ambient manifold. Substituting this into the evolution equation above yields 
\begin{align*} 
\frac{d}{dt} \bigg ( \int_{M_t} f_{\sigma,+}^p \bigg ) 
&\leq -\frac{1}{C} \, p(p-1) \int_{M_t} f_{\sigma,+}^{p-2} \, |\nabla f_\sigma|^2 - \frac{1}{C} \, p \int_{M_t} f_{\sigma,+}^p \, \frac{|\nabla h|^2}{G_\kappa^2} \\ 
&+ C\sigma p^2 \int_{M_t} f_{\sigma,+}^{p-1} \, \frac{|\nabla h|}{G_\kappa} \, |\nabla f_\sigma| + Cp \int_{M_t} G_\kappa^\sigma \, f_{\sigma,+}^{p-1} - \frac{1}{C} \int_{M_t} G_\kappa^2 \, f_{\sigma,+}^p, 
\end{align*} 
provided that $p$ is sufficiently large and $\sigma$ is sufficiently small. This implies 
\begin{align*} 
\frac{d}{dt} \bigg ( \int_{M_t} f_{\sigma,+}^p \bigg ) 
&\leq -\frac{1}{C} \, p(p-1) \int_{M_t} f_{\sigma,+}^{p-2} \, |\nabla f_\sigma|^2 - \frac{1}{C} \, p \int_{M_t} f_{\sigma,+}^p \, \frac{|\nabla h|^2}{G_\kappa^2} \\ 
&+ Cp \int_{M_t} G_\kappa^\sigma \, f_{\sigma,+}^{p-1} - \frac{1}{C} \int_{M_t} G_\kappa^2 \, f_{\sigma,+}^p, 
\end{align*} 
provided that $p$ is sufficiently large and $\sigma p^{\frac{1}{2}}$ is sufficiently small. Since $G_\kappa$ is uniformly bounded from below on bounded time intervals, we have 
\[Cp \, G_\kappa^\sigma \, f_{\sigma,+}^{p-1} - \frac{1}{C} \, G_\kappa^2 \, f_{\sigma,+}^p \leq (C' p)^p \, G_\kappa^{2-(2-\sigma) \, p} \leq (C'' p)^p.\] 
This completes the proof of Proposition \ref{Lp.bound}. \\ 

\begin{corollary} 
\label{cor.of.Lp.bound}
Assume that $p \geq \frac{1}{c_0}$ and $\sigma \leq c_0 \, p^{-\frac{1}{2}}$. Then we have 
\[\int_{M_t} f_{\sigma,+}^p \leq C,\] 
where $C$ is a positive constant that depends only on $p$, $\sigma$, $\delta$, $T$, $\kappa$, the initial hypersurface $M_0$, and the ambient manifold.
\end{corollary}

We now continue with the proof of Theorem \ref{cylindrical.estimate}. For $k \geq 0$, we define 
\[f_{\sigma,k} = G_\kappa^{\sigma-1} \, \Big ( H - \frac{(n-1)^2(n+2)}{4} \, (1+\delta) \, G_\kappa \Big ) - k\] 
and 
\[f_{\sigma,k,+} = \max \{f_{\sigma,k},0\}.\] 

\begin{proposition}
\label{final.ingredient}
We have
\begin{align*}
\frac{d}{dt} \bigg ( \int_{M_t} f_{\sigma,k,+}^p \bigg )
&\leq -\frac{1}{C} \, p(p-1) \int_{M_t} f_{\sigma,k,+}^{p-2} \, |\nabla f_{\sigma,k}|^2 \\
&+ C \sigma \, p \int_{M_t} G_\kappa^2 \, f_{\sigma,k,+}^{p-1} \, f_\sigma + (Cp)^p \, |M_t \cap \{f_{\sigma,k} \geq 0\}|
\end{align*}
if $k \geq 0$ and $p \geq \frac{1}{c_1}$. Here, $c_1$ and $C$ are a positive constants that depend only on $\delta$, $T$, $\kappa$, the initial hypersurface $M_0$, and the ambient manifold.
\end{proposition}

\textbf{Proof.}
Assume that $k \geq 0$. The function $f_{\sigma,k}$ satisfies 
\begin{align*}
\frac{\partial}{\partial t} f_{\sigma,k}
&\leq \frac{\partial G_\kappa}{\partial h_{ij}} \, D_i D_j f_{\sigma,k} + C \, \frac{|\nabla G_\kappa|}{G_\kappa} \, |\nabla f_{\sigma,k}| \\ 
&+ C \, \sigma \, G_\kappa^2 \, f_\sigma - \frac{1}{C} \, G_\kappa^{\sigma-2} \, |\nabla h|^2 + C \, G_\kappa^\sigma. 
\end{align*} 
This implies 
\begin{align*}
&\frac{d}{dt} \bigg ( \int_{M_t} f_{\sigma,k,+}^p \bigg ) \\
&\leq p \int_{M_t} f_{\sigma,k,+}^{p-1} \, \frac{\partial G_\kappa}{\partial h_{ij}} \, D_i D_j f_{\sigma,k} + Cp \int_{M_t} f_{\sigma,k,+}^{p-1} \, \frac{|\nabla G_\kappa|}{G_\kappa} \, |\nabla f_{\sigma,k}| \\ 
&+ C \sigma p \int_{M_t} G_\kappa^2 \, f_{\sigma,k,+}^{p-1} \, f_\sigma - \frac{1}{C} \, p \int_{M_t} G_\kappa^{\sigma-2} \, f_{\sigma,k,+}^{p-1} \, |\nabla h|^2 \\ 
&+ Cp \int_{M_t} G_\kappa^\sigma \, f_{\sigma,k,+}^{p-1} - \frac{1}{C} \int_{M_t} G_\kappa^2 \, f_{\sigma,k,+}^p.  
\end{align*} 
As above, integration by parts yields 
\begin{align*} 
&\frac{d}{dt} \bigg ( \int_{M_t} f_{\sigma,k,+}^p \bigg ) \\ 
&\leq -p(p-1) \int_{M_t} f_{\sigma,k,+}^{p-2} \, \frac{\partial G_\kappa}{\partial h_{ij}} \, D_i f_{\sigma,k} \, D_j f_{\sigma,k} + Cp \int_{M_t} f_{\sigma,k,+}^{p-1} \, \frac{|\nabla h|}{G_\kappa} \, |\nabla f_{\sigma,k}| \\ 
&+ C \sigma p \int_{M_t} G_\kappa^2 \, f_{\sigma,k,+}^{p-1} \, f_\sigma - \frac{1}{C} \, p \int_{M_t} f_{\sigma,k,+}^p \, \frac{|\nabla h|^2}{G_\kappa^2} \\ 
&+ Cp \int_{M_t} G_\kappa^\sigma \, f_{\sigma,k,+}^{p-1} - \frac{1}{C} \int_{M_t} G_\kappa^2 \, f_{\sigma,k,+}^p, 
\end{align*} 
hence 
\begin{align*} 
&\frac{d}{dt} \bigg ( \int_{M_t} f_{\sigma,k,+}^p \bigg ) \\ 
&\leq -\frac{1}{C} \, p(p-1) \int_{M_t} f_{\sigma,k,+}^{p-2} \, \frac{\partial G_\kappa}{\partial h_{ij}} \, D_i f_{\sigma,k} \, D_j f_{\sigma,k} \\ 
&+ C \sigma p \int_{M_t} G_\kappa^2 \, f_{\sigma,k,+}^{p-1} \, f_\sigma + Cp \int_{M_t} G_\kappa^\sigma \, f_{\sigma,k,+}^{p-1} - \frac{1}{C} \int_{M_t} G_\kappa^2 \, f_{\sigma,k,+}^p 
\end{align*} 
for $p$ sufficiently large. Finally, we have 
\[Cp \, G_\kappa^\sigma \, f_{\sigma,+}^{p-1} - \frac{1}{C} \, G_\kappa^2 \, f_{\sigma,+}^p \leq (C' p)^p \, G_\kappa^{2-(2-\sigma) \, p} \leq (C'' p)^p.\] From this, the assertion follows. \\

We now complete the proof of Theorem \ref{cylindrical.estimate}. To that end, we show that $f_\sigma$ is uniformly bounded from above for some small number $\sigma>0$. The proof uses Stampacchia iteration. Let us fix real numbers $p$ and $\sigma$ such that $p \geq \frac{1}{\min\{c_0,c_1\}}$ and $0 < \sigma < c_0 \, (2np)^{-\frac{1}{2}} - 2 \, p^{-1}$. For abbreviation, let $A(k) = \int_0^T |M_t \cap \{f_{\sigma,k} \geq 0\}|$. It follows from Proposition \ref{final.ingredient} that 
\begin{align*}
\frac{d}{dt} \bigg ( \int_{M_t} f_{\sigma,k,+}^p \bigg )
&\leq -\frac{1}{C} \int_{M_t} f_{\sigma,k,+}^{p-2} \, |\nabla f_{\sigma,k}|^2 \\
&+ C \int_{M_t} G_\kappa^2 \, f_{\sigma,k,+}^{p-1} \, f_\sigma + C \, |M_t \cap \{f_{\sigma,k} \geq 0\}|, 
\end{align*} 
where $C$ is a positive constant that depends only on $p$, $\sigma$, $\delta$, $T$, $\kappa$, $M_0$, and $X$, but not on $k$. If $k \geq K_0 := \max \{ \sup_{M_0} \frac{H}{G_\kappa},\sup_{M_0} H\}$, then we have $f_{\sigma,k} \leq 0$ on the initial hypersurface $M_0$. This implies 
\[\sup_{t \in [0,T)} \int_{M_t} f_{\sigma,k,+}^p \leq C \, A(k) + C \int_0^T \int_{M_t} G_\kappa^2 \, f_{\sigma,k,+}^{p-1} \, f_\sigma\] 
and 
\[\int_0^T \int_{M_t} f_{\sigma,k,+}^{p-2} \, |\nabla f_{\sigma,k}|^2 \leq C \, A(k) + C \int_0^T \int_{M_t} G_\kappa^2 \, f_{\sigma,k,+}^{p-1} \, f_\sigma\]
for $k \geq K_0$. Here, $C$ is a positive constant that depends only on $p$, $\sigma$, $\delta$, $T$, $\kappa$, $M_0$, and $X$, but on $k$. Using the Michael-Simon Sobolev inequality (cf. \cite{Michael-Simon}), we obtain 
\begin{align*} 
\bigg ( \int_{M_t} f_{\sigma,k,+}^{\frac{pn}{n-1}} \bigg )^{\frac{n-1}{n}} 
&\leq C \int_{M_t} f_{\sigma,k,+}^{p-1} \, |\nabla f_{\sigma,k}| + C \int_{M_t} (G_\kappa+1) \, f_{\sigma,k,+}^p \\ 
&\leq C \int_{M_t} f_{\sigma,k,+}^{p-2} \, |\nabla f_{\sigma,k}|^2 + C \int_{M_t} (G_\kappa^2+1) \, f_{\sigma,k,+}^p 
\end{align*} 
for $k \geq K_0$. Integrating over $t$ gives 
\begin{align*} 
\int_0^T \bigg ( \int_{M_t} f_{\sigma,k,+}^{\frac{pn}{n-1}} \bigg )^{\frac{n-1}{n}} 
&\leq C \int_0^T \int_{M_t} f_{\sigma,k,+}^{p-2} \, |\nabla f_{\sigma,k}|^2 + C \int_0^T \int_{M_t} (G_\kappa^2+1) \, f_{\sigma,k,+}^p \\ 
&\leq C \, A(k) + C \int_0^T \int_{M_t} (G_\kappa^2+1) \, f_{\sigma,k,+}^{p-1} \, f_\sigma 
\end{align*} 
for $k \geq K_0$. Hence, it follows from H\"older's inequality that 
\begin{align*} 
&\bigg ( \int_0^T \int_{M_t} f_{\sigma,k,+}^{\frac{p(n+1)}{n}} \bigg )^{\frac{n}{n+1}} \\ 
&\leq \bigg ( \sup_{t \in [0,T)} \int_{M_t} f_{\sigma,k,+}^p \bigg )^{\frac{1}{n+1}} \cdot  \bigg ( \int_0^T \bigg ( \int_{M_t} f_{\sigma,k,+}^{\frac{pn}{n-1}} \bigg )^{\frac{n-1}{n}} \bigg )^{\frac{n}{n+1}} \\ 
&\leq C \, A(k) + C \int_0^T \int_{M_t} (G_\kappa^2+1) \, f_{\sigma,k,+}^{p-1} \, f_\sigma
\end{align*} 
for $k \geq K_0$. As above, $C$ is a positive constant that depends only on $p$, $\sigma$, $\delta$, $T$, $\kappa$, $M_0$, and $X$, but not on $k$. By Corollary \ref{cor.of.Lp.bound}, we have 
\[\int_0^T \int_{M_t} (G_\kappa^{4n}+1) \, f_{\sigma,+}^{2np} \leq C \int_0^T \int_{M_t} (f_{\sigma+2p^{-1},+}^{2np} + f_{\sigma,+}^{2np}) \leq C.\] 
Applying H\"older's inequality again, we obtain 
\begin{align*} 
&\bigg ( \int_0^T \int_{M_t} f_{\sigma,k,+}^{\frac{p(n+1)}{n}} \bigg )^{\frac{n}{n+1}} \\
&\leq C \, A(k) + C \int_0^T \int_{M_t} (G_\kappa^2+1) \, f_{\sigma,+}^p \, 1_{\{f_{\sigma,k} \geq 0\}} \\ 
&\leq C \, A(k) + C \, \bigg ( \int_0^T \int_{M_t} (G_\kappa^{4n}+1) \, f_{\sigma,+}^{2np} \bigg )^{\frac{1}{2n}} \, \bigg ( \int_0^T \int_{M_t} 1_{\{f_{\sigma,k} \geq 0\}} \bigg )^{1-\frac{1}{2n}} \\ 
&\leq C \, A(k)^{1-\frac{1}{2n}} 
\end{align*} 
for $k \geq K_0$. Thus, we conclude that  
\[A(\tilde{k})^{1-\frac{1}{n+1}} \, (\tilde{k}-k)^p \leq C \, A(k)^{1-\frac{1}{2n}}\] 
for $\tilde{k} \geq k \geq K_0$. Again, $C$ is a positive constant that depends only on $p$, $\sigma$, $\delta$, $T$, $\kappa$, $M_0$, and $X$, but not on $k$ or $\tilde{k}$. Iterating this inequality gives $A(k)=0$ for some constant $k = k(p,\sigma,\delta,T,\kappa,M_0,X)$. From this, we deduce that
\[H \leq \frac{(n-1)^2(n+2)}{4} \, (1+2\delta) \, G_\kappa + B,\] 
where $B$ is a positive constant that depends only on $\delta$, $T$, $\kappa$, the initial hypersurface $M_0$, and the ambient manifold. This completes the proof of Theorem \ref{cylindrical.estimate}. \\

\begin{proposition}
\label{alg}
We have 
\[\frac{3(n-2)}{n+2} \, \lambda_1 \geq \frac{(n-1)^2(n+2)}{4} \, G_\kappa - H + \frac{(n-1)(n+6)}{n+2} \, \kappa.\] 
\end{proposition} 

\textbf{Proof.} 
We define 
\[a_{ij} = \begin{cases} 1 & \text{\rm if $1<i<j$} \\ 2 & \text{\rm if $1=i<j$.} \end{cases}\] 
Using the Cauchy-Schwarz inequality, we obtain 
\begin{align*} 
\frac{(n-1)^2(n+2)^2}{4} 
&= \Big ( \sum_{i<j} a_{ij} \Big )^2 \\ 
&\leq \Big ( \sum_{i<j} \frac{1}{\lambda_i+\lambda_j-2\kappa} \Big ) \, \Big ( \sum_{i<j} a_{ij}^2 \, (\lambda_i+\lambda_j-2\kappa) \Big )^2 \\ 
&= G_\kappa^{-1} \, \Big ( \sum_{1<i<j} (\lambda_i+\lambda_j-2\kappa) + 4 \, \sum_{1<j} (\lambda_1+\lambda_j-2\kappa) \Big ) \\ 
&= G_\kappa^{-1} \, ((n+2) \, H + 3(n-2) \, \lambda_1 - (n-1)(n+6)\kappa). 
\end{align*} 
This proves the assertion. \\ 

By combining Theorem \ref{cylindrical.estimate} with Proposition \ref{alg}, we obtain an analogue of the convexity estimates for mean curvature flow established by the second author and Carlo Sinestrari \cite{Huisken-Sinestrari1},\cite{Huisken-Sinestrari2} (see also \cite{Andrews-Langford-McCoy2}, where a different class of fully nonlinear flows is studied). 

\begin{corollary}[Convexity Estimate] 
\label{convexity.estimate}
Let $M_t$, $t \in [0,T)$, be a family of closed, $\kappa$-two-convex hypersurfaces moving with speed $G_\kappa$, and let $\delta$ be an arbitrary positive real number. Then 
\[\lambda_1 \geq -\delta \, G_\kappa - C,\] 
where $C$ is a positive constant that depends only on $\delta$, $T$, $\kappa$, the initial hypersurface $M_0$, and the ambient manifold.
\end{corollary}

The following result is similar in spirit to Hamilton's strict maximum principle for the Ricci flow (cf. \cite{Hamilton1}).

\begin{proposition}[Splitting Theorem]
\label{splitting}
Suppose that $M_t$, $t \in [-1,0]$, is a family of (possibly non-closed) two-convex hypersurfaces in $\mathbb{R}^{n+1}$ which move with velocity $G = \big ( \sum_{i<j} \frac{1}{\lambda_i+\lambda_j} \big )^{-1}$. Moreover, suppose that $M_t$ satisfies the pointwise inequality $H \leq \frac{(n-1)^2(n+2)}{4} \, G$. Then either $\lambda_1>0$ at each point in the interior of $M_0$, or else each hypersurface $M_t$ is contained in a cylinder.
\end{proposition}

\textbf{Proof.}
Suppose that $\lambda_1 \leq 0$ at some point in the interior of $M_0$. At that point, we have $H \geq \frac{(n-1)^2(n+2)}{4} \, G$ by Proposition \ref{alg}. Using the strict maximum principle, we conclude that $H = \frac{(n-1)^2(n+2)}{4} \, G$ at all points in spacetime. On the other hand, we have 
\[\frac{\partial}{\partial t} G = \frac{\partial G}{\partial h_{ij}} \, (D_i D_j G + h_{ik} \, h_{jk} \, G)\] 
and 
\[\frac{\partial}{\partial t} H \leq \frac{\partial G}{\partial h_{ij}} \, (D_i D_j H + h_{ik} \, h_{jk} \, H) - \frac{1}{C} \, \frac{|\nabla h|^2}{G}\] 
in view of Lemma \ref{good.term}. Since $H$ is a constant multiple of $G$, we conclude that $|\nabla h|^2 = 0$ at each point on in spacetime. In other words, the second fundamental form is parallel. Therefore, $M_t$ is contained in a cylinder. \\

\section{The inscribed radius estimate}

\label{inscr.rad}

Let $M_t$, $t \in [0,T)$, be a family of embedded hypersurfaces in a compact Riemannian manifold which move with velocity $G_\kappa$. For each point on $M_t$, the inscribed radius is defined as the radius of the largest geodesic ball which is contained in $\Omega_t$ and touches $M_t$ at that point. 

It will be convenient to parametrise the hypersurfaces $M_t$ by a map $F: M \times [0,T) \to X$. We define 
\[\mu(x,t) = \sup_{y \in M, \, 0 < d(F(x,t),F(y,t)) \leq \frac{1}{2} \, \text{\rm inj}(X)} \Big ( -\frac{2 \, \langle \exp_{F(x,t)}^{-1}(F(y,t)),\nu(x,t) \rangle}{d(F(x,t),F(y,t))^2} \Big ).\] 
For hypersurfaces in Euclidean space, $\mu$ is equal to the reciprocal of the inscribed radius at the point $(x,t)$. When $\kappa=0$ and $X=\mathbb{R}^{n+1}$, Andrews, Langford, and McCoy \cite{Andrews-Langford-McCoy1} established an important estimate for the inscribed radius along the flow. Their work directly generalizes to the case $\kappa \geq 0$. The estimate can also be extended to the Riemannian setting:

\begin{proposition}
\label{evol}
Consider a point $(\bar{x},\bar{t}) \in M \times [0,T)$ such that $\lambda_n(\bar{x},\bar{t}) < \mu(\bar{x},\bar{t})$ and $\mu(\bar{x},\bar{t})$ is sufficiently large. We further assume that $\Phi: M \times [0,\bar{t}] \to \mathbb{R}$ is a smooth function such that $\Phi(\bar{x},\bar{t}) = \mu(\bar{x},\bar{t})$ and $\Phi(x,t) \geq \mu(x,t)$ for all points $(x,t) \in M \times [0,\bar{t}]$. Then 
\[\frac{\partial \Phi}{\partial t} \leq \sum_{i,j} \frac{\partial G_\kappa}{\partial h_{ij}} \, (D_i D_j \Phi + h_{ik} \, h_{jk} \, \Phi) + C \, \Phi + C \, \sum_i \frac{1}{\Phi-\lambda_i}\] 
at the point $(\bar{x},\bar{t})$. Here, $C$ is a positive constant that depends only on $T$, $\kappa$, the initial hypersurface $M_0$, and the ambient manifold.
\end{proposition}

\textbf{Proof.} 
We sketch the details for the convenience of the reader. For each point $q \in X$, we define a function $\psi_q: X \to \mathbb{R}$ by $\psi_q(p) = \frac{1}{2} \, d(p,q)^2$, where $d(p,q)$ denotes the Riemannian distance in $X$. For abbreviation, we put $\Xi_{q,p} := (\text{\rm Hess} \, \psi_q)_p - g$. Clearly, $\Xi_{q,p}$ is a symmetric bilinear form on $T_p X$, and we have $|\Xi_{q,p}| \leq O(d(p,q)^2)$. We define a function $Z: M \times M \times [0,T) \to \mathbb{R}$ by 
\begin{align*} 
Z(x,y,t) 
&= \Phi(x,t) \, \psi_{F(y,t)}(F(x,t)) - \big \langle \nabla \psi_{F(y,t)} \big |_{F(x,t)},\nu(x,t) \big \rangle \\ 
&= \frac{1}{2} \, \Phi(x,t) \, d(F(x,t),F(y,t))^2 + \big \langle \exp_{F(x,t)}^{-1}(F(y,t)),\nu(x,t) \big \rangle. 
\end{align*} 
By assumption, we have $Z(x,y,t) \geq 0$ whenever $x \in M$, $t \in [0,\bar{t}]$, and $d(F(x,t),F(y,t)) \leq \frac{1}{2} \, \text{\rm inj}(X)$. Moreover, we can find a point $\bar{y} \in M$ such that $0 < d(F(\bar{x},\bar{t}),F(\bar{y},\bar{t})) \leq \frac{1}{2} \, \text{\rm inj}(X)$ and $Z(\bar{x},\bar{y},\bar{t}) = 0$. Clearly, 
\[\Phi(\bar{x},\bar{t}) \, d(F(\bar{x},\bar{t}),F(\bar{y},\bar{t})) \leq 2.\] 
Let us choose geodesic normal coordinates around $\bar{x}$ such that $h_{ij}(\bar{x},\bar{t})$ is a diagonal matrix. Moreover, we put $\lambda_i = h_{ii}(\bar{x},\bar{t})$ and $\gamma_i = \frac{\partial G_\kappa}{\partial \lambda_i}$. The first variation of $Z$ with respect to $x$ gives  
\begin{align*}
0 = \frac{\partial Z}{\partial x_i}(\bar{x},\bar{y},\bar{t}) 
&= \frac{1}{2} \, \frac{\partial \Phi}{\partial x_i}(\bar{x},\bar{t}) \, d(F(\bar{x},\bar{t}),F(\bar{y},\bar{t}))^2 \\ 
&- \Phi(\bar{x},\bar{t}) \, \Big \langle \exp_{F(\bar{x},\bar{t})}^{-1}(F(\bar{y},\bar{t})),\frac{\partial F}{\partial x_i}(\bar{x},\bar{t}) \Big \rangle \\ 
&+ h_i^j(\bar{x},\bar{t}) \, \Big \langle \exp_{F(\bar{x},\bar{t})}^{-1}( F(\bar{y},\bar{t})),\frac{\partial F}{\partial x_j}(\bar{x},\bar{t}) \Big \rangle \\ 
&- \Xi_{F(\bar{y},\bar{t}),F(\bar{x},\bar{t})} \Big ( \frac{\partial F}{\partial x_i}(\bar{x},\bar{t}),\nu(\bar{x},\bar{t}) \Big ). 
\end{align*} 
Consequently,  
\begin{align*} 
&\Big \langle \exp_{F(\bar{x},\bar{t})}^{-1}( F(\bar{y},\bar{t})),\frac{\partial F}{\partial x_i}(\bar{x},\bar{t}) \Big \rangle \\ 
&= \frac{1}{2} \, \frac{1}{\Phi(\bar{x},\bar{t}) - \lambda_i(\bar{x},\bar{t})} \, \Big ( \frac{\partial \Phi}{\partial x_i}(\bar{x},\bar{t}) + O(1) \Big ) \, d(F(\bar{x},\bar{t}),F(\bar{y},\bar{t}))^2. 
\end{align*} 
We next consider the second variation of $Z$ with respect to $x$. Using the Codazzi equations, we obtain 
\begin{align*}
&\sum_i \gamma_i \, \frac{\partial^2 Z}{\partial x_i^2}(\bar{x},\bar{y},\bar{t}) \\ 
&= \frac{1}{2} \sum_i \gamma_i \, \frac{\partial^2 \Phi}{\partial x_i^2}(\bar{x},\bar{t}) \, d(F(\bar{x},\bar{t}),F(\bar{y},\bar{t}))^2 \\ 
&- 2 \sum_i \gamma_i \, \frac{\partial \Phi}{\partial x_i}(\bar{x},\bar{t}) \, \Big \langle \exp_{F(\bar{x},\bar{t})}^{-1}( F(\bar{y},\bar{t})),\frac{\partial F}{\partial x_i}(\bar{x},\bar{t}) \Big \rangle \\ 
&+ \sum_i \frac{\partial G_\kappa}{\partial x_i}(\bar{x},\bar{t}) \, \Big \langle \exp_{F(\bar{x},\bar{t})}^{-1}( F(\bar{y},\bar{t})),\frac{\partial F}{\partial x_i}(\bar{x},\bar{t}) \Big \rangle \\ 
&+ \sum_i \gamma_i \, \lambda_i \, \Phi(\bar{x},\bar{t}) \, \langle \exp_{F(\bar{x},\bar{t})}^{-1}( F(\bar{y},\bar{t})),\nu(\bar{x},\bar{t}) \rangle \\
&- \sum_i \gamma_i \, \lambda_i^2 \, \langle \exp_{F(\bar{x},\bar{t})}^{-1}( F(\bar{y},\bar{t})),\nu(\bar{x},\bar{t}) \rangle \\  
&+ \Phi(\bar{x},\bar{t}) \, \sum_i \gamma_i - \sum_i \gamma_i \, \lambda_i \\ 
&+ O \big ( d(F(\bar{x},\bar{t}),F(\bar{y},\bar{t})) \big ), 
\end{align*} 
hence 
\begin{align*}
&\sum_i \gamma_i \, \frac{\partial^2 Z}{\partial x_i^2}(\bar{x},\bar{y},\bar{t}) \\ 
&\leq \frac{1}{2} \, \bigg ( \sum_i \gamma_i \, \frac{\partial^2 \Phi}{\partial x_i^2}(\bar{x},\bar{t}) + \sum_i \gamma_i \, \lambda_i^2 \, \Phi(\bar{x},\bar{t}) \\ 
&\hspace{20mm} - \sum_i \frac{2}{\Phi(\bar{x},\bar{t})-\lambda_i(\bar{x},\bar{t})} \, \gamma_i \, \Big ( \frac{\partial \Phi}{\partial x_i}(\bar{x},\bar{t}) \Big )^2 \bigg ) \, d(F(\bar{x},\bar{t}),F(\bar{y},\bar{t}))^2 \\ 
&+ \sum_i \frac{\partial G_\kappa}{\partial x_i}(\bar{x},\bar{t}) \, \Big \langle \exp_{F(\bar{x},\bar{t})}^{-1}(F(\bar{y},\bar{t})),\frac{\partial F}{\partial x_i}(\bar{x},\bar{t}) \Big \rangle \\ 
&+ \sum_i \gamma_i \, \lambda_i \, \Phi(\bar{x},\bar{t}) \, \big \langle \exp_{F(\bar{x},\bar{t})}^{-1}(F(\bar{y},\bar{t})),\nu(\bar{x},\bar{t}) \big \rangle \\
&+ \Phi(\bar{x},\bar{t}) \, \sum_i \gamma_i - \sum_i \gamma_i \, \lambda_i \\ 
&+ O \bigg ( d(F(\bar{x},\bar{t}),F(\bar{y},\bar{t})) + \sum_i \frac{1}{\Phi(\bar{x},\bar{t})-\lambda_i(\bar{x},\bar{t})} \, \Big | \frac{\partial \Phi}{\partial x_i}(\bar{x},\bar{t}) \Big | \, d(F(\bar{x},\bar{t}),F(\bar{y},\bar{t}))^2 \bigg ). 
\end{align*} 
For a suitable choice of the coordinate system around $\bar{y}$, we have
\[\frac{\partial^2 Z}{\partial x_i \, \partial y_i}(\bar{x},\bar{y},\bar{t}) = -(\Phi(\bar{x},\bar{t})-\lambda_i(\bar{x},\bar{t})) + O \big ( d(F(\bar{x},\bar{t}),F(\bar{y},\bar{t})) \big )\] 
for each $i$. Hence, for this choice of the coordinate system around $\bar{y}$, we obtain  
\[\sum_i \gamma_i \, \frac{\partial^2 Z}{\partial x_i \, \partial y_i}(\bar{x},\bar{y},\bar{t}) = -\Phi(\bar{x},\bar{t}) \, \sum_i \gamma_i + \sum_i \gamma_i \, \lambda_i + O \big ( d(F(\bar{x},\bar{t}),F(\bar{y},\bar{t})) \big ).\] 
Finally, the second variation of $Z$ with respect to $y$ is given by 
\begin{align*} 
&\frac{\partial^2 Z}{\partial y_i^2}(\bar{x},\bar{y},\bar{t}) \\ 
&= \Phi(\bar{x},\bar{t}) - h_{ii}(\bar{y},\bar{t}) \, \big \langle (D\exp_{F(\bar{x},\bar{t})}^{-1})_{F(\bar{y},\bar{t})}(\nu(\bar{y},\bar{t})),\nu(\bar{x},\bar{t}) + \Phi(\bar{x},\bar{t}) \, \exp_{F(\bar{x},\bar{t})}^{-1}(F(\bar{y},\bar{t})) \big \rangle \\ 
&+ O \big ( d(F(\bar{x},\bar{t}),F(\bar{y},\bar{t})) \big ). 
\end{align*} 
Note that 
\begin{align*} 
&\big \langle (D\exp_{F(\bar{x},\bar{t})}^{-1})_{F(\bar{y},\bar{t})}(\nu(\bar{y},\bar{t})),\nu(\bar{x},\bar{t}) + \Phi(\bar{x},\bar{t}) \, \exp_{F(\bar{x},\bar{t})}^{-1}(F(\bar{y},\bar{t})) \big \rangle \\ 
&= 1 + O \big ( d(F(\bar{x},\bar{t}),F(\bar{y},\bar{t}))^2 \big ). 
\end{align*}
As in Andrews-Langford-McCoy \cite{Andrews-Langford-McCoy1}, we have
\[G_\kappa(\bar{y},\bar{t}) \leq G_\kappa(\bar{x},\bar{t}) + \sum_i \gamma_i \, (h_{ii}(\bar{y},\bar{t})-\lambda_i)\] 
since $G_\kappa$ is concave. This gives 
\begin{align*} 
&\sum_i \gamma_i \, \frac{\partial^2 Z}{\partial y_i^2}(\bar{x},\bar{y},\bar{t}) \\ 
&\leq \Phi(\bar{x},\bar{t}) \, \sum_i \gamma_i - \sum_i \gamma_i \, \lambda_i + G_\kappa(\bar{x},\bar{t}) \\ 
&- G_\kappa(\bar{y},\bar{t}) \, \big \langle (D\exp_{F(\bar{x},\bar{t})}^{-1})_{F(\bar{y},\bar{t})}(\nu(\bar{y},\bar{t})),\nu(\bar{x},\bar{t}) + \Phi(\bar{x},\bar{t}) \, \exp_{F(\bar{x},\bar{t})}^{-1}(F(\bar{y},\bar{t})) \big \rangle \\ 
&+ O \big ( d(F(\bar{x},\bar{t}),F(\bar{y},\bar{t})) \big ). 
\end{align*} 
Putting these facts together yields 
\begin{align*}
&\sum_i \gamma_i \, \Big ( \frac{\partial^2 Z}{\partial x_i^2}(\bar{x},\bar{y},\bar{t}) + 2 \, \frac{\partial^2 Z}{\partial x_i \, \partial y_i}(\bar{x},\bar{y},\bar{t}) + \frac{\partial^2 Z}{\partial y_i^2}(\bar{x},\bar{y},\bar{t}) \Big ) \\ 
&\leq \frac{1}{2} \, \bigg ( \sum_i \gamma_i \, \frac{\partial^2 \Phi}{\partial x_i^2}(\bar{x},\bar{t}) + \sum_i \gamma_i \, \lambda_i^2 \, \Phi(\bar{x},\bar{t}) \\ 
&\hspace{20mm} - \sum_i \frac{2}{\Phi(\bar{x},\bar{t})-\lambda_i(\bar{x},\bar{t})} \, \gamma_i \, \Big ( \frac{\partial \Phi}{\partial x_i}(\bar{x},\bar{t}) \Big )^2 \bigg ) \, d(F(\bar{x},\bar{t}),F(\bar{y},\bar{t}))^2 \\ 
&+ \sum_i \frac{\partial G_\kappa}{\partial x_i}(\bar{x},\bar{t}) \, \Big \langle \exp_{F(\bar{x},\bar{t})}^{-1}(F(\bar{y},\bar{t})),\frac{\partial F}{\partial x_i}(\bar{x},\bar{t}) \Big \rangle \\ 
&+ G_\kappa(\bar{x},\bar{t}) + \sum_i \gamma_i \, \lambda_i \,  \Phi(\bar{x},\bar{t}) \, \big \langle \exp_{F(\bar{x},\bar{t})}^{-1}(F(\bar{y},\bar{t})),\nu(\bar{x},\bar{t}) \big \rangle \\ 
&- G_\kappa(\bar{y},\bar{t}) \, \big \langle (D\exp_{F(\bar{x},\bar{t})}^{-1})_{F(\bar{y},\bar{t})}(\nu(\bar{y},\bar{t})),\nu(\bar{x},\bar{t}) + \Phi(\bar{x},\bar{t}) \, \exp_{F(\bar{x},\bar{t})}^{-1}(F(\bar{y},\bar{t})) \big \rangle \\ 
&+ O \bigg ( d(F(\bar{x},\bar{t}),F(\bar{y},\bar{t})) + \sum_i \frac{1}{\Phi(\bar{x},\bar{t})-\lambda_i(\bar{x},\bar{t})} \, \Big | \frac{\partial \Phi}{\partial x_i}(\bar{x},\bar{t}) \Big | \, d(F(\bar{x},\bar{t}),F(\bar{y},\bar{t}))^2 \bigg ). 
\end{align*} 
On the other hand, we have 
\begin{align*} 
\frac{\partial Z}{\partial t}(\bar{x},\bar{y},\bar{t}) 
&= \frac{1}{2} \, \frac{\partial \Phi}{\partial t}(\bar{x},\bar{t}) \, d(F(\bar{x},\bar{t}),F(\bar{y},\bar{t}))^2 \\ 
&+ G_\kappa(\bar{x},\bar{t}) + G_\kappa(\bar{x},\bar{t}) \, \Phi(\bar{x},\bar{t}) \, \big \langle \exp_{F(\bar{x},\bar{t})}^{-1}(F(\bar{y},\bar{t})),\nu(\bar{x},\bar{t}) \big \rangle \\
&- G_\kappa(\bar{y},\bar{t}) \, \big \langle (D\exp_{F(\bar{x},\bar{t})}^{-1})_{F(\bar{y},\bar{t})}(\nu(\bar{y},\bar{t})),\nu(\bar{x},\bar{t}) + \Phi(\bar{x},\bar{t}) \, \exp_{F(\bar{x},\bar{t})}^{-1}(F(\bar{y},\bar{t})) \big \rangle \\ 
&+ \sum_i \frac{\partial G_\kappa}{\partial x_i}(\bar{x},\bar{t}) \, \Big \langle \exp_{F(\bar{x},\bar{t})}^{-1}(F(\bar{y},\bar{t})),\frac{\partial F}{\partial x_i}(\bar{x},\bar{t}) \Big \rangle \\ 
&+ G_\kappa(\bar{x},\bar{t}) \, \Xi_{F(\bar{y},\bar{t}),F(\bar{x},\bar{t})}(\nu(\bar{x},\bar{t}),\nu(\bar{x},\bar{t})). 
\end{align*} 
Finally, we have $G_\kappa(\bar{x},\bar{t}) \leq \sum_i \gamma_i \, \lambda_i$ and $\big \langle \exp_{F(\bar{x},\bar{t})}^{-1}(F(\bar{y},\bar{t})),\nu(\bar{x},\bar{t}) \big \rangle \leq 0$. This implies 
\begin{align*}
0 &\geq \frac{\partial Z}{\partial t}(\bar{x},\bar{y},\bar{t}) - \sum_i \gamma_i \, \Big ( \frac{\partial^2 Z}{\partial x_i^2}(\bar{x},\bar{y},\bar{t}) + 2 \, \frac{\partial^2 Z}{\partial x_i \, \partial y_i}(\bar{x},\bar{y},\bar{t}) + \frac{\partial^2 Z}{\partial y_i^2}(\bar{x},\bar{y},\bar{t}) \Big ) \\ 
&\geq \frac{1}{2} \, \bigg ( \frac{\partial \Phi}{\partial t}(\bar{x},\bar{t}) - \sum_i \gamma_i \, \frac{\partial^2 \Phi}{\partial x_i^2}(\bar{x},\bar{t}) - \sum_i \gamma_i \, \lambda_i^2 \, \Phi(\bar{x},\bar{t}) \\ 
&\hspace{20mm} + \sum_i \frac{2}{\Phi(\bar{x},\bar{t})-\lambda_i(\bar{x},\bar{t})} \, \gamma_i \, \Big ( \frac{\partial \Phi}{\partial x_i}(\bar{x},\bar{t}) \Big )^2 \bigg ) \, d(F(\bar{x},\bar{t}),F(\bar{y},\bar{t}))^2 \\ 
&- O \bigg ( d(F(\bar{x},\bar{t}),F(\bar{y},\bar{t})) + \sum_i \frac{1}{\Phi(\bar{x},\bar{t})-\lambda_i(\bar{x},\bar{t})} \, \Big | \frac{\partial \Phi}{\partial x_i}(\bar{x},\bar{t}) \Big | \, d(F(\bar{x},\bar{t}),F(\bar{y},\bar{t}))^2 \bigg ). 
\end{align*} 
Finally, we multiply both sides of the previous inequality by $\frac{2}{d(F(\bar{x},\bar{t}),F(\bar{y},\bar{t}))^2}$. Using the estimate 
\begin{align*} 
&\frac{1}{d(F(\bar{x},\bar{t}),F(\bar{y},\bar{t}))} \\ 
&\leq \frac{|\langle \exp_{F(\bar{x},\bar{t})}^{-1}(F(\bar{y},\bar{t})),\nu(\bar{x},\bar{t}) \rangle|}{d(F(\bar{x},\bar{t}),F(\bar{y},\bar{t}))^2} + \sum_i \frac{|\langle \exp_{F(\bar{x},\bar{t})}^{-1}(F(\bar{y},\bar{t})),\frac{\partial F}{\partial x_i}(\bar{x},\bar{t}) \rangle|}{d(F(\bar{x},\bar{t}),F(\bar{y},\bar{t}))^2} \\ 
&\leq \frac{1}{2} \, \Phi(\bar{x},\bar{t}) + \sum_i \frac{1}{2} \, \frac{1}{\Phi(\bar{x},\bar{t}) - \lambda_i(\bar{x},\bar{t})} \, \Big ( \Big | \frac{\partial \Phi}{\partial x_i}(\bar{x},\bar{t}) \Big | + O(1) \Big ), 
\end{align*} 
we obtain  
\begin{align*}
&\frac{\partial \Phi}{\partial t}(\bar{x},\bar{t}) - \sum_i \gamma_i \, \frac{\partial^2 \Phi}{\partial x_i^2}(\bar{x},\bar{t}) - \sum_i \gamma_i \, \lambda_i^2 \, \Phi(\bar{x},\bar{t}) \\ 
&+ \sum_i \frac{2}{\Phi(\bar{x},\bar{t})-\lambda_i(\bar{x},\bar{t})} \, \gamma_i \, \Big ( \frac{\partial \Phi}{\partial x_i}(\bar{x},\bar{t}) \Big )^2 \\ 
&\leq O \bigg ( \Phi(\bar{x},\bar{t}) + \sum_i \frac{1}{\Phi(\bar{x},\bar{t}) - \lambda_i(\bar{x},\bar{t})} + \sum_i \frac{1}{\Phi(\bar{x},\bar{t}) - \lambda_i(\bar{x},\bar{t})} \, \Big | \frac{\partial \Phi}{\partial x_i}(\bar{x},\bar{t}) \Big | \bigg ). 
\end{align*} From this, the assertion follows. \\

\begin{corollary} 
\label{noncollapsing.riemannian.case}
At each point on $M_t$, the inscribed radius is bounded from below by $\frac{\alpha}{G_\kappa}$, where $\alpha$ is a positive constant that depends only on $T$, $\kappa$, the initial hypersurface $M_0$, and the ambient manifold.
\end{corollary} 

\textbf{Proof.} 
By Proposition \ref{evol}, the function $\mu$ satisfies 
\[\frac{\partial}{\partial t} \mu \leq \sum_{i,j} \frac{\partial G_\kappa}{\partial h_{ij}} \, (D_i D_j \mu + h_{ik} \, h_{jk} \, \mu) + C \, \mu + C \, \sum_i \frac{1}{\mu-\lambda_i}\] 
whenever $\mu>\lambda_n$ and $\mu$ is sufficiently large. Here, the inequality is interpreted in the viscosity sense. Furthermore, $C$ is a positive constant that depends only on $T$, $\kappa$, the initial hypersurface $M_0$, and the ambient manifold. We next observe that $\lambda_i \leq H \leq \beta^{-1} \, G_\kappa$ by Proposition \ref{uniformly.two.convex}. This implies 
\[\frac{\partial}{\partial t} \mu \leq \sum_{i,j} \frac{\partial G_\kappa}{\partial h_{ij}} \, (D_i D_j \mu + h_{ik} \, h_{jk} \, \mu) + C \, \mu\] 
whenever $\frac{\mu}{G_\kappa}$ is sufficiently large. By the maximum principle, the ratio $\frac{\mu}{G_\kappa}$ is uniformly bounded from above on bounded time intervals. Since the inscribed radius is comparable to $\frac{1}{\mu}$, the assertion follows. \\

\section{Interior estimates for radial graphs}

\label{rad.graph}

In this section, we establish interior estimates for solutions of the fully nonlinear flow which can be written as radial graphs. These estimates are similar in spirit to the interior estimates for graphs evolving by mean curvature flow proved by Klaus Ecker and the second author \cite{Ecker-Huisken1},\cite{Ecker-Huisken2}; see also \cite{Gerhardt}, \cite{Tso}, and \cite{Urbas}, where global estimates for radial graphs evolving under other fully nonlinear curvature flows are established.

\begin{proposition} 
\label{interior.curvature.estimate.for.radial.graphs}
Let $X$ be a Riemannian manifold, let $p$ be a point on $X$, and let $r \leq \min \{1,\frac{1}{4} \, \text{\rm inj}(X)\}$. Suppose that $\Omega_t$, $t \in [-r^2,0]$, is a one-parameter family of smooth open domains in $X$ such that $B_r(p) \subset \Omega_t$ and the hypersurfaces $M_t = \partial \Omega_t$ move with velocity $G_\kappa$. Finally, we assume that $\langle -\exp_x^{-1}(p),\nu \rangle \geq 10^{-3} \, r$ and $G_\kappa \geq \beta H$ at each point $x \in \partial U_t \cap B_{2r}(p)$, where $U_t$ denotes the connected component of $\Omega_t \cap B_{2r}(p)$ which contains the ball $B_r(p)$. Then the norm of the second fundamental form satisfies 
\[\Big ( \frac{25r^2}{9}-d(p,x)^2 \Big ) \, (t+r^2)^{\frac{1}{2}} \, |h| \leq C \, r^2\] 
for all $t \in [-r^2,0]$ and all $x \in \partial U_t \cap B_{\frac{5r}{3}}(p)$. Here, $C$ is a positive constant that depends only on $\beta$ and the ambient manifold.
\end{proposition}

Of course, the number $10^{-3}$ in the statement of Proposition \ref{interior.curvature.estimate.for.radial.graphs} can be replaced by any positive constant. \\

\textbf{Proof.} 
The function $G_\kappa$ satisfies 
\[\frac{\partial}{\partial t} G_\kappa \leq \frac{\partial G_\kappa}{\partial h_{ij}} \, (D_i D_j G_\kappa + h_{ik} \, h_{jk} \, G_\kappa) + C \, G_\kappa.\] 
We next consider the radial vector field $V = -\exp_x^{-1}(p)$ on the ball $B_{2r}(p)$. The function $\langle V,\nu \rangle$ satisfies 
\[\frac{\partial}{\partial t} \langle V,\nu \rangle \geq \frac{\partial G_\kappa}{\partial h_{ij}} \, (D_i D_j \langle V,\nu \rangle + h_{ik} \, h_{jk} \, \langle V,\nu \rangle) - C \, |h| - C.\] 
We next define $v = (\langle V,\nu \rangle^2 - \sigma r^2)^{-\frac{1}{2}}$, where $\sigma = 10^{-7}$. By assumption, the product $rv$ is uniformly bounded from above and below for each $t \in [-r^2,0]$ and each $x \in \partial U_t \cap B_{2r}(p)$. Moreover, the function $v$ satisfies 
\begin{align*} 
\frac{\partial}{\partial t} v 
&\leq \frac{\partial G_\kappa}{\partial h_{ij}} \, (D_i D_j v - h_{ik} \, h_{jk} \, (v+\sigma r^2 v^3)) \\ 
&- \Big ( 3 - \frac{1}{1+\sigma r^2 v^2} \Big ) \, v^{-1} \, \frac{\partial G_\kappa}{\partial h_{ij}} \, D_i v \, D_j v \\ 
&+ C \, r^{-2} \, |h| + C \, r^{-2} 
\end{align*} 
for each $t \in [-r^2,0]$ and each $x \in \partial U_t \cap B_{2r}(p)$. Finally, the function $\eta = \frac{25r^2}{9}-d(p,x)^2$ satisfies 
\[\frac{\partial}{\partial t} \eta \leq \frac{\partial G_\kappa}{\partial h_{ij}} \, D_i D_j \eta + C.\] 
Hence, if we put $\psi = \eta \, v \, G_\kappa$, then we obtain 
\begin{align*} 
\frac{\partial}{\partial t} \psi 
&\leq \frac{\partial G_\kappa}{\partial h_{ij}} \, D_i D_j \psi - 2 \, \eta \, \frac{\partial G_\kappa}{\partial h_{ij}} \, D_i v \, D_j G_\kappa - 2 \, v \, \frac{\partial G_\kappa}{\partial h_{ij}} \, D_i \eta \, D_j G_\kappa \\ 
&- 2 \, G_\kappa \, \frac{\partial G_\kappa}{\partial h_{ij}} \, D_i \eta \, D_j v - \eta \, G_\kappa \, \Big ( 3 - \frac{1}{1+\sigma r^2v^2} \Big ) \, v^{-1} \, \frac{\partial G_\kappa}{\partial h_{ij}} \, D_i v \, D_j v \\ 
&- \sigma r^2 \, \eta \, v^3 \, G_\kappa \, \frac{\partial G_\kappa}{\partial h_{ij}} \, h_{ik} \, h_{jk} + C \, |h| \, G_\kappa + C \, r^{-1} \, G_\kappa \\ 
&= \frac{\partial G_\kappa}{\partial h_{ij}} \, D_i D_j \psi - 2 \, v^{-1} \, \frac{\partial G_\kappa}{\partial h_{ij}} \, D_i v \, D_j \psi - 2 \, \eta^{-1} \, \frac{\partial G_\kappa}{\partial h_{ij}} \, D_i \eta \, D_j \psi \\ 
&+ 2 \, G_\kappa \, \frac{\partial G_\kappa}{\partial h_{ij}} \, D_i \eta \, D_j v + 2 \, v \, G_\kappa \, \eta^{-1} \, \frac{\partial G_\kappa}{\partial h_{ij}} \, D_i \eta \, D_j \eta - \eta \, G_\kappa \, \frac{\sigma r^2 v}{1+\sigma r^2v^2} \, \frac{\partial G_\kappa}{\partial h_{ij}} \, D_i v \, D_j v \\ 
&- \sigma r^2 \, \eta \, v^3 \, G_\kappa \, \frac{\partial G_\kappa}{\partial h_{ij}} \, h_{ik} \, h_{jk} + C \, |h| \, G_\kappa + C \, r^{-1} \, G_\kappa \\ 
&= \frac{\partial G_\kappa}{\partial h_{ij}} \, D_i D_j \psi - 2 \, v^{-1} \, \frac{\partial G_\kappa}{\partial h_{ij}} \, D_i v \, D_j \psi - 2 \, \eta^{-1} \, \frac{\partial G_\kappa}{\partial h_{ij}} \, D_i \eta \, D_j \psi \\ 
&- \frac{\sigma r^2 \eta v}{1+\sigma r^2 v^2} \, G_\kappa \, \frac{\partial G_\kappa}{\partial h_{ij}} \, \Big ( D_i v - \frac{1+\sigma r^2 v^2}{\sigma r^2 \eta v} \, D_i \eta \Big ) \, \Big ( D_j v - \frac{1+\sigma r^2 v^2}{\sigma r^2 \eta v} \, D_j \eta \Big ) \\ 
&+ \frac{1+\sigma r^2 v^2}{\sigma r^2 v} \, G_\kappa \, \eta^{-1} \, \frac{\partial G_\kappa}{\partial h_{ij}} \, D_i \eta \, D_j \eta + 2 \, v \, G_\kappa \, \eta^{-1} \, \frac{\partial G_\kappa}{\partial h_{ij}} \, D_i \eta \, D_j \eta \\ 
&- \sigma r^2 \, \eta \, v^3 \, G_\kappa \, \frac{\partial G_\kappa}{\partial h_{ij}} \, h_{ik} \, h_{jk} + C \, |h| \, G_\kappa + C \, r^{-1} \, G_\kappa 
\end{align*} 
for each $t \in [-r^2,0]$ and each $x \in \partial U_t \cap B_{\frac{5r}{3}}(p)$. Clearly, $\frac{\partial G_\kappa}{\partial h_{ij}} \, D_i \eta \, D_j \eta \leq C(n) \, r^2$. Moreover, our assumptions imply that $|h| \leq C \, G_\kappa$. Furthermore, it follows from the Cauchy-Schwarz inequality that $G_\kappa^2 \leq \big ( \sum_{i,j} \frac{\partial G_\kappa}{\partial h_{ij}} \, h_{ij} \big )^2 \leq C(n) \sum_{i,j,k} \frac{\partial G_\kappa}{\partial h_{ij}} \, h_{ik} \, h_{jk}$. Putting these facts together, we obtain 
\begin{align*} 
\frac{\partial}{\partial t} \psi 
&\leq \frac{\partial G_\kappa}{\partial h_{ij}} \, D_i D_j \psi - 2 \, v^{-1} \, \frac{\partial G_\kappa}{\partial h_{ij}} \, D_i v \, D_j \psi - 2 \, \eta^{-1} \, \frac{\partial G_\kappa}{\partial h_{ij}} \, D_i \eta \, D_j \psi \\ 
&- \frac{1}{C} \, r^{-1} \, \eta \, G_\kappa^3 + C \, G_\kappa^2 + C \, r \, \eta^{-1} \, G_\kappa, 
\end{align*} 
hence 
\begin{align*} 
\frac{\partial}{\partial t} \psi 
&\leq \frac{\partial G_\kappa}{\partial h_{ij}} \, D_i D_j \psi - 2 \, v^{-1} \, \frac{\partial G_\kappa}{\partial h_{ij}} \, D_i v \, D_j \psi - 2 \, \eta^{-1} \, \frac{\partial G_\kappa}{\partial h_{ij}} \, D_i \eta \, D_j \psi \\ 
&- \frac{1}{C} \, r^2 \, \eta^{-2} \, \psi^3 + C \, r^2 \, \eta^{-2} \, \psi^2 + C \, r^2 \, \eta^{-2} \, \psi 
\end{align*} 
for each $t \in [-r^2,0]$ and each $x \in \partial U_t \cap B_{\frac{5r}{3}}(p)$. We now define 
\[Q(t) = \sup_{x \in \partial U_t \cap B_{\frac{5r}{3}}(p)} \psi(x,t)\] 
for $t \in [-r^2,0]$. If $Q(t)$ is sufficiently large, then we have 
\[\frac{\partial}{\partial t} \psi \leq -\frac{1}{C} \, r^2 \, \eta^{-2} \, \psi^3 \leq -\frac{1}{C} \, r^{-2} \, \psi^3\] 
for each point $x \in \partial U_t \cap B_{\frac{5r}{3}}(p)$ satisfying $\psi(x,t) = Q(t)$. Hence, if $Q(t)$ is sufficiently large, then we obtain 
\[\limsup_{t' \nearrow t} \frac{Q(t)-Q(t')}{t-t'} \leq -\frac{1}{C} \, r^{-2} \, Q(t)^3.\] 
This finally gives 
\[Q(t) \leq C \, r \, (t+r^2)^{-\frac{1}{2}}\] 
for all $t \in [-r^2,0]$. Since $|h| \leq C \, G_\kappa$, the assertion follows. \\

\begin{corollary} 
\label{higher.regularity.for.radial.graphs}
Let $X$ be a Riemannian manifold, let $p$ be a point in $X$, and let $r \leq \min \{1,\frac{1}{4} \, \text{\rm inj}(X)\}$. Suppose that $\Omega_t$, $t \in [-r^2,0]$, is a one-parameter family of smooth open domains in $X$ such that $B_r(p) \subset \Omega_t$ and the hypersurfaces $M_t = \partial \Omega_t$ move with velocity $G_\kappa$. Finally, we assume that $\langle -\exp_x^{-1}(p),\nu \rangle \geq 10^{-3} \, r$ and $G_\kappa \geq \beta H$ at each point $x \in \partial U_t \cap B_{2r}(p)$, where $U_t$ denotes the connected component of $\Omega_t \cap B_{2r}(p)$ which contains the ball $B_r(p)$. Then 
\[r^2 \, |\nabla h(x,0)| + r^3 \, |\nabla^2 h(x,0)| \leq \Lambda\] 
for all points $x \in \partial U_0 \cap B_{\frac{4r}{3}}(p)$ satisfying $G_\kappa(x,0) \geq \alpha \, r^{-1}$. Here, $\Lambda$ is a positive constant that depends only on $\alpha$, $\beta$, and the ambient manifold.
\end{corollary}

\textbf{Proof.} 
By Proposition \ref{interior.curvature.estimate.for.radial.graphs}, we can find a positive constant $K \geq 100$ such that $|h| \leq K \, r^{-1}$ and $G_\kappa \leq K \, r^{-1}$ for all $t \in [-\frac{r^2}{4},0]$ and all $x \in \partial U_t \cap B_{\frac{3r}{2}}(p)$. We now fix a point $x \in \partial U_0 \cap B_{\frac{4r}{3}}(p)$ satisfying $G_\kappa(x,0) \geq \alpha \, r^{-1}$. For each $t \in [-\frac{r^2}{100K},0]$, we have $\mathcal{P}(x,0,\frac{r}{100},\frac{r^2}{100K}) \cap M_t \subset \partial U_t \cap B_{\frac{3r}{2}}(p)$. In particular, we have $|h| \leq K \, r^{-1}$ and $G_\kappa \leq K \, r^{-1}$ at each point in $\mathcal{P}(x,0,\frac{r}{100},\frac{r^2}{100K})$. Hence, on the set $\mathcal{P}(x,0,\frac{r}{100},\frac{r^2}{100K})$, the function $G_\kappa$ satisfies a uniformly parabolic equation with bounded coefficients. 

Using the Krylov-Safonov theorem (cf. Theorem \ref{krylov.safonov.interior.estimate}), we obtain a H\"older estimate for $G_\kappa$ on the set $\mathcal{P}(x,0,\frac{r}{200},\frac{r^2}{200K})$. In particular, there exists a uniform constant $\theta \in (0,\frac{1}{400K})$ such that $\frac{\alpha}{2} \, r^{-1} \leq G_\kappa \leq K \, r^{-1}$ on the set $\mathcal{P}(x,0,\theta \, r,\theta \, r^2)$. Theorem \ref{interior.estimate.for.fully.nonlinear.PDE} now gives H\"older estimates for the second fundamental form on the set $\mathcal{P}(x,0,\frac{\theta}{2} \, r,\frac{\theta}{2} \, r^2)$. Using Schauder theory, we obtain estimates for all derivatives of the second fundamental form on the set $\mathcal{P}(x,0,\frac{\theta}{4} \, r,\frac{\theta}{4} \, r^2)$. In particular, this gives bounds for $|\nabla h(x,0)|$ and $|\nabla^2 h(x,0)|$. \\

\section{The pointwise curvature derivative estimate}

\label{derivative.est}

In this section, we establish a pointwise estimate for the derivatives of the second fundamental form. We begin by introducing some notation. Let $\varphi(s) = \tan(\frac{1}{100}) \, (s+s^2)$ for $s \in [0,1]$. Given two points $p$ and $x$ satisfying $d(p,x) < \frac{1}{2} \, \text{\rm inj}(X)$, we define 
\begin{align*} 
C_{p,x} 
&= \{\exp_x(s \exp_x^{-1}(p) + v): s \in (0,1), \, v \in T_x X, \, \langle \exp_x^{-1}(p),v \rangle = 0, \\ 
&\hspace{45mm} |v| < \varphi(s) \, d(p,x)\} 
\end{align*} 
and 
\begin{align*} 
S_{p,x} 
&= \{\exp_x(s \exp_x^{-1}(p) + v): s \in (0,1), \, v \in T_x X, \, \langle \exp_x^{-1}(p),v \rangle = 0, \\ 
&\hspace{45mm} |v| = \varphi(s) \, d(p,x)\}. 
\end{align*} 
It is easy to see that 
\[\partial C_{p,x} \subset S_{p,x} \cup \{x\} \cup B_{\frac{1}{4} d(p,x)}(p)\] 
whenever $d(p,x)$ is sufficiently small. Near $x$, the hypersurface $S_{p,x}$ is asymptotic to a cone with aperture $2 \cdot \frac{1}{100}$. Note that $S_{p,x}$ is slightly bent outwards, as a consequence of the choice of the function $\varphi(s)$. We refer to $C_{p,x}$ as a pseudo-cone.

\begin{lemma} 
\label{smallest.curvature.eigenvalue.of.S}
If $d(p,x)$ is sufficiently small, then, at each point on $S_{p,x}$, the smallest curvature eigenvalue is less than $-10^{-3} \, d(p,x)^{-1}$.
\end{lemma} 

\textbf{Proof.} 
The smallest curvature eigenvalue of $S_{p,x}$ is given by 
\[-(1+\varphi'(s))^{-\frac{3}{2}} \, \varphi''(s) \, d(p,x)^{-1} + O(d(p,x)).\] 
Since $(1+\varphi'(s))^{-\frac{3}{2}} \, \varphi''(s) > 10^{-3}$ for all $s \in [0,1]$, the assertion follows. \\

Suppose now that $\Omega_t$, $t \in [0,T)$, is a one-parameter family of smooth open domains in $X$ with the property that the hypersurfaces $M_t = \partial \Omega_t$ move with velocity $G_\kappa$. By Proposition \ref{uniformly.two.convex}, we have $G_\kappa \geq \beta H$, where $\beta$ is a positive constant that depends only on $T$, $\kappa$, the initial hypersurface $M_0$, and the ambient manifold. Moreover, by Corollary \ref{noncollapsing.riemannian.case}, there exists a constant $\alpha > 0$, depending only on $T$, $\kappa$, the initial hypersurface $M_0$, and the ambient manifold, such that the inscribed radius is at least $\frac{\alpha}{G_\kappa}$ at each point in spacetime. 

The following is the main result of this section:

\begin{theorem} 
\label{curvature.derivative.estimate}
We have $\alpha^2 \, G_\kappa^{-2} \, |\nabla h| + \alpha^3 \, G_\kappa^{-3} \, |\nabla^2 h(x,0)| \leq \Lambda$ whenever $G_\kappa$ is sufficiently large. Here, $\alpha$ is the constant in Corollary \ref{noncollapsing.riemannian.case}, and $\Lambda$ is the constant appearing in Corollary \ref{higher.regularity.for.radial.graphs}.
\end{theorem} 

\textbf{Proof.} 
Suppose that the assertion is false. Then there exists a sequence of points $(x_k,t_k)$ in spacetime such that $G_\kappa(x_k,t_k) \to \infty$ and 
\[\alpha^2 \, G_\kappa(x_k,t_k)^{-2} \, |\nabla h(x_k,t_k)| + \alpha^3 \, G_\kappa(x_k,t_k)^{-3} \, |\nabla^2 h(x_k,t_k)| > \Lambda\] 
for each $k$. Using a standard point-picking argument, we can find, for each $k$, a point $(\bar{x}_k,\bar{t}_k)$ with the following properties: 
\begin{itemize}
\item[(i)] $\bar{t}_k \leq t_k$. 
\item[(ii)] $G_\kappa(\bar{x}_k,\bar{t}_k) \geq G_\kappa(x_k,t_k)$.
\item[(iii)] $\alpha^2 \, G_\kappa(\bar{x}_k,\bar{t}_k)^{-2} \, |\nabla h(\bar{x}_k,\bar{t}_k)| + \alpha^3 \, G_\kappa(\bar{x}_k,\bar{t}_k)^{-3} \, |\nabla^2 h(\bar{x}_k,\bar{t}_k)| > \Lambda$.
\item[(iv)] $\alpha^2 \, G_\kappa(x,t)^{-2} \, |\nabla h(x,t)| + \alpha^3 \, G_\kappa(x,t)^{-3} \, |\nabla^2 h(x,t)| \leq \Lambda$ for all points $(x,t)$ with $t \leq \bar{t}_k$ and $G_\kappa(x,t) \geq 2 \, G_\kappa(\bar{x}_k,\bar{t}_k)$. 
\end{itemize} 
For abbreviation, let $r_k = \alpha \, G_\kappa(\bar{x}_k,\bar{t}_k)^{-1}$. Note that $r_k \to 0$ in view of property (ii). Using Corollary \ref{noncollapsing.riemannian.case}, we can find a point $p_k$ such that $\bar{x}_k \in \partial B_{r_k}(p_k)$ and $B_{r_k}(p_k) \subset \Omega_{\bar{t}_k}$. Clearly, $B_{r_k}(p_k) \subset \Omega_t$ for all $t \in [\bar{t}_k-r_k^2,\bar{t}_k]$. For each $t \in [\bar{t}_k-r_k^2,\bar{t}_k]$, we denote by $U_t^{(k)}$ the connected component of $\Omega_t \cap B_{2r_k}(p_k)$ which contains the ball $B_{r_k}(p_k)$. Clearly, the sets $U_t^{(k)}$ shrink as $t$ increases. We distinguish two cases:

\textit{Case 1:} Suppose that $C_{p_k,x} \subset U_t^{(k)}$ for all $t \in [\bar{t}_k-r_k^2,\bar{t}_k]$ and all points $x \in U_t^{(k)}$. This implies $\varangle (-\exp_x^{-1}(p_k),\nu) \leq \frac{\pi}{2}-\frac{1}{100}$ for all $t \in [\bar{t}_k-r_k^2,\bar{t}_k]$ and all points $x \in \partial U_t^{(k)} \cap B_{2r_k}(p_k)$. Consequently, $\langle -\exp_x^{-1}(p_k),\nu \rangle \geq 10^{-3} \, r_k$ for all $t \in [\bar{t}_k-r_k^2,\bar{t}_k]$ and all points $x \in \partial U_t^{(k)} \cap B_{2r_k}(p_k)$. Using Corollary \ref{higher.regularity.for.radial.graphs}, we obtain $r_k^2 \, |\nabla h(x,\bar{t}_k)| + r_k^3 \, |\nabla^2 h(x,\bar{t}_k)| \leq \Lambda$ for all points $x \in \partial U_{\bar{t}_k}^{(k)} \cap B_{\frac{4r_k}{3}}(p_k)$. On the other hand, we clearly have $\bar{x}_k \in \partial U_{\bar{t}_k}^{(k)} \cap B_{\frac{4r_k}{3}}(p_k)$, $G_\kappa(\bar{x}_k,\bar{t}_k) = \alpha \, r_k^{-1}$, and furthermore $r_k^2 \, |\nabla h(\bar{x}_k,\bar{t}_k)| + r_k^3 \, |\nabla^2 h(\bar{x}_k,\bar{t}_k)| > \Lambda$ in view of property (iii) above. This is a contradiction.

\textit{Case 2:} Suppose that there exists a time $\tilde{t}_k \in [\bar{t}_k-r_k^2,\bar{t}_k]$ and a point $x \in U_{\tilde{t}_k}^{(k)}$ such that $C_{p_k,x} \not\subset U_{\tilde{t}_k}^{(k)}$. Let 
\[A^{(k)} = \{x \in U_{\tilde{t}_k}^{(k)}: C_{p_k,x} \subset U_{\tilde{t}_k}^{(k)}\}.\] 
It is clear that $B_{r_k}(p_k) \subset A^{(k)}$, $A^{(k)} \neq U_{\tilde{t}_k}^{(k)}$, and $A^{(k)}$ is relatively closed as a subset of $U_{\tilde{t}_k}^{(k)}$. Since $U_{\tilde{t}_k}^{(k)}$ is connected, it follows that $A^{(k)}$ cannot be an open set. Consequently, there exists a point $\tilde{x}_k \in A^{(k)}$ with the property that $B_\sigma(\tilde{x}_k) \not\subset A^{(k)}$ for all $\sigma>0$. Note that $r_k \leq d(p_k,\tilde{x}_k) < 2r_k$ and $C_{p_k,\tilde{x}_k} \subset U_{\tilde{t}_k}^{(k)}$. Moreover, since 
\[\partial C_{p_k,\tilde{x}_k} \subset S_{p_k,\tilde{x}_k} \cup \{\tilde{x}_k\} \cup B_{\frac{r_k}{2}}(p_k) \subset S_{p_k,\tilde{x}_k} \cup U_{\tilde{t}_k}^{(k)},\] 
the hypersurface $S_{p_k,\tilde{x}_k}$ touches $\partial U_{\tilde{t}_k}^{(k)} \cap B_{2r_k}(p_k)$ somewhere from the inside. 

Let us consider a point $y_k$ where the hypersurface $S_{p_k,\tilde{x}_k}$ touches $\partial U_{\tilde{t}_k}^{(k)} \cap B_{2r_k}(p_k)$ from the inside. Clearly, $r_k \leq d(p_k,y_k) < 2r_k$. Since $C_{p_k,\tilde{x}_k}$ has aperture $2 \cdot \frac{1}{100}$, we can find a unit vector $v_k \in T_{y_k} X$ such that $\varangle (v_k,\nu(y_k,\tilde{t}_k)) \geq \frac{\pi}{2}+10^{-3}$ and 
\[\{\exp_{y_k}(sv_k): 0 < s < \frac{r_k}{2}\} \subset C_{p_k,\tilde{x}_k} \subset U_{\tilde{t}_k}^{(k)}.\] 
In particular, we have 
\begin{equation} 
\label{a}
\{\exp_{y_k}(sv_k): 0 < s < \frac{r_k}{2}\} \cap M_{\tilde{t}_k} = \emptyset. 
\end{equation}
By Lemma \ref{smallest.curvature.eigenvalue.of.S}, the smallest curvature eigenvalue of $S_{p_k,\tilde{x}_k}$ is less than $-10^{-3} \, d(p_k,\tilde{x}_k)^{-1}$ at each point on $S_{p_k,\tilde{x}_k}$. Since the hypersurface $S_{p_k,\tilde{x}_k}$ touches $\partial U_{\tilde{t}_k}^{(k)} \cap B_{2r_k}(p_k)$ from the inside at $y_k$, it follows that 
\[\lambda_1(y_k,\tilde{t}_k) \leq -10^{-3} \, d(p_k,\tilde{x}_k)^{-1} \leq -\frac{10^{-3}}{2} \, r_k^{-1} = -\frac{10^{-3}}{2} \, \alpha^{-1} \, G_\kappa(\bar{x}_k,\bar{t}_k).\] 
In particular, $\lambda_1(y_k,\tilde{t}_k) \to -\infty$ in view of property (ii) above. Using Corollary \ref{convexity.estimate}, we obtain 
\[\lambda_1(y_k,\tilde{t}_k) \, G_\kappa(y_k,\tilde{t}_k)^{-1} \to 0.\] 
Thus, we conclude that 
\[G_\kappa(y_k,\tilde{t}_k) \, G_\kappa(\bar{x}_k,\bar{t}_k)^{-1} \to \infty.\] 
In particular, we have $G_\kappa(y_k,\tilde{t}_k) \geq 8 \, G_\kappa(\bar{x}_k,\bar{t}_k)$ if $k$ is sufficiently large. For each $k$, we define 
\[L_k = \min \Big \{ \inf \Big \{ G_\kappa(y_k,\tilde{t}_k) \, d_{M_{\tilde{t}_k}}(y_k,x): x \in M_{\tilde{t}_k}, \, \frac{G_\kappa(x,\tilde{t}_k)}{G_\kappa(y_k,\tilde{t}_k)} \notin [\frac{1}{2},2] \Big \},10^6 \Big \}.\] 
By definition of $L_k$, we have $\frac{1}{2} \, G_\kappa(y_k,\tilde{t}_k) \leq G_\kappa(x,\tilde{t}_k) \leq 2 \, G_\kappa(y_k,\tilde{t}_k)$ for all points $x \in M_{\tilde{t}_k}$ satisfying $d_{M_{\tilde{t}_k}}(y_k,x) \leq L_k \, G_\kappa(y_k,\tilde{t}_k)^{-1}$. Using property (iv) above, we obtain 
\[\sup_{\mathcal{P}(y_k,\tilde{t}_k,(L_k+\theta) \, G_\kappa(y_k,\tilde{t}_k)^{-1},\theta \, G_\kappa(y_k,\tilde{t}_k)^{-2})} G_\kappa \leq 4 \, G_\kappa(y_k,\tilde{t}_k)\] 
and 
\[\inf_{\mathcal{P}(y_k,\tilde{t}_k,(L_k+\theta) \, G_\kappa(y_k,\tilde{t}_k)^{-1},\theta \, G_\kappa(y_k,\tilde{t}_k)^{-2})} G_\kappa \geq \frac{1}{4} \, G_\kappa(y_k,\tilde{t}_k),\] 
where $\theta$ is a positive constant independent of $k$.

In the next step, we restrict the flow to the parabolic neighborhood $\mathcal{P}(y_k,\tilde{t}_k,(L_k+\theta) \, G_\kappa(y_k,\tilde{t}_k)^{-1},\theta \, G_\kappa(y_k,\tilde{t}_k)^{-2})$. On this parabolic neighborhood, the ratio $\frac{H}{G_\kappa(y_k,\tilde{t}_k)}$ is uniformly bounded from above, and the ratio $\frac{\lambda_1+\lambda_2-2\kappa}{G_\kappa(y_k,\tilde{t}_k)}$ is uniformly bounded from below. Hence, if we perform a parabolic dilation around the point $(y_k,\tilde{t}_k)$ with factor $G_\kappa(y_k,\tilde{t}_k)$, then the rescaled flow has bounded curvature and is uniformly two-convex. By property (iv) above, we have bounds for the first and second derivatives of the second fundamental form. Hence, the rescaled flows converge to a smooth, non-flat limit flow in $\mathbb{R}^{n+1}$, which moves with normal velocity $G$ and satisfies the pointwise inequality $H \leq \frac{(n-1)^2(n+2)}{4} \, G$ (see Theorem \ref{cylindrical.estimate}). Since $\lambda_1(y_k,\tilde{t}_k) < 0$ for each $k$, there exists a point on the limit flow where the smallest curvature eigenvalue is non-positive. Using Proposition \ref{splitting}, we conclude that the limit flow is contained in a family of shrinking cylinders. In particular, this implies 
\[\sup \{G_\kappa(x,\tilde{t}_k): x \in M_{\tilde{t}_k}, \, d_{M_{\tilde{t}_k}}(y_k,x) \leq L_k \, G_\kappa(y_k,\tilde{t}_k)^{-1}\} \leq (1+o(1)) \, G_\kappa(y_k,\tilde{t}_k)\] 
and 
\[\inf \{G_\kappa(x,\tilde{t}_k): x \in M_{\tilde{t}_k}, \, d_{M_{\tilde{t}_k}}(y_k,x) \leq L_k \, G_\kappa(y_k,\tilde{t}_k)^{-1}\} \geq (1-o(1)) \, G_\kappa(y_k,\tilde{t}_k).\] 
Consequently, we have $L_k=10^6$ if $k$ is sufficiently large. Moreover, the point $y_k$ lies at the center of an $(\varepsilon_k,6,\frac{(n-1)(n+2)}{4} \cdot 10^5)$-neck in $M_{\tilde{t}_k}$ for some sequence $\varepsilon_k \to 0$. Since $\varangle (v_k,\nu(y_k,\tilde{t}_k)) \geq \frac{\pi}{2}+10^{-3}$, we conclude that 
\begin{equation} 
\label{b}
\{\exp_{y_k}(sv_k): 0 < s < 10^4 \, G_\kappa(y_k,\tilde{t}_k)^{-1}\} \cap M_{\tilde{t}_k} \neq \emptyset 
\end{equation} 
if $k$ is sufficiently large. Since $G_\kappa(y_k,\tilde{t}_k) \, r_k \to \infty$, the statements (\ref{a}) and (\ref{b}) are in contradiction. This completes the proof of Theorem \ref{curvature.derivative.estimate}. \\

\section{A-priori estimates for surgically modified flows}

\label{estimates.for.flows.with.surgery}

In this section, we consider flows with velocity $G_\kappa$ which are interrupted by finitely many surgeries. We first explain how some basic notions introduced in \cite{Huisken-Sinestrari2} can be adapted to the Riemannian setting. 

\begin{definition} 
Suppose that $M$ is a hypersurface in a Riemannian manifold $X$, and let $p$ be a point in $M$. We say that $p$ lies at the center of an $(\varepsilon,k,L)$-neck in $M$ if $0$ lies at the center of an $(\varepsilon,k,L)$-neck in $\exp_p^{-1}(M \cap B_{\frac{1}{4} \text{\rm inj}(X)}(p)) \subset T_p X$ in the sense of Defintion 3.1 (v) in \cite{Huisken-Sinestrari2}. 
\end{definition}

By a result of Hamilton \cite{Hamilton4}, a neck admits a canonical foliation by spheres which have constant mean curvature with respect to the induced metric on the neck. If the radius of the neck is sufficiently small, each leaf of Hamilton's foliation bounds a unique area-minimizing disk in ambient space; this gives a canonical foliation of the solid tube associated with the neck (see \cite{Huisken-Sinestrari2}, Proposition 3.25). We next define the axis of the neck. There are several ways of doing this. For example, to each leaf $\Sigma$ in Hamilton's foliation, we may associate a point $z \in X$ such that set $\exp_z^{-1}(\Sigma) \subset T_z X$ has its center of mass at the origin. The collection of all these points corresponding to different leaves of the foliation is a smooth curve, which we call the axis of the neck. 

We next explain how to do surgery on a neck, when the ambient space is a Riemannian manifold. As before, let $M$ be a hypersurface in a Riemannian manifold $X$, and suppose that $N \subset M$ is a neck in $M$. To perform surgery on such a neck, we pick a point $z$ on the axis of $N$. It is easy to see that $\exp_z^{-1}(N) \subset T_z X$ is a neck in Euclidean space. On this neck, we may perform a standard surgery as defined on pp.~154-155 in \cite{Huisken-Sinestrari2}. As a result, we obtain a capped-off neck in $T_z X$. We then paste the capped-off neck back into $X$ using the exponential map $\exp_z$.

With this understood, we can now give a precise definition of a surgically modified flow.

\begin{definition} 
A surgically modified flow is a family of closed, embedded, $\kappa$-two-convex hypersurfaces $M_t = \partial \Omega_t$, $t \in [0,T)$, with the following properties: 
\begin{itemize}
\item The hypersurfaces $M_t$ move smoothly with speed $G_\kappa = \big ( \sum_{i<j} \frac{1}{\lambda_i+\lambda_j-2\kappa} \big )^{-1}$, except at finitely many times.
\item At each of these times, we perform finitely many standard surgeries. Each surgery is performed in the middle third of an $(\varepsilon,6,\frac{(n-1)(n+2)}{4} \, L)$-neck, where $L \geq 10^9$. On each neck on which surgery is being performed, the curvature satisfies $\frac{1}{2} \, G_* \leq G_\kappa \leq 2G_*$, where $G_*$ is a large positive number (the same for all surgeries).
\item During each surgery procedure, we glue in a cap. The construction of this cap is described in detail in \cite{Huisken-Sinestrari2}. In particular, the intrinsic diameter of the cap is less than $100 \, G_*^{-1}$. Moreover, we have $\frac{1}{2} \, G_* \leq G_\kappa \leq 100 \, G_*$ and $G_*^{-1} \, |h| + G_*^{-2} \, |\nabla h| + G_*^{-3} \, |\nabla^2 h| + G_*^{-4} \, |\nabla^3 h| \leq C(n)$ at each point on the cap. 
\item Immediately after surgery, some components may be removed. Each of these components bounds a region which is diffeomorphic to $B^n$ or $B^{n-1} \times S^1$.
\end{itemize} 
The number $G_*$ will be referred to as the surgery scale of the flow $M_t$.
\end{definition}

\begin{lemma}
\label{properties.of.surgery}
We can find surgery parameters $B,\tau_0$ and positive numbers $G_\Diamond,\varepsilon_\Diamond,\sigma_\Diamond$ such that the following holds. Suppose that we perform a standard surgery with parameters $B,\tau_0$ on an $(\varepsilon,6,\frac{(n-1)(n+2)}{4} \, L)$-neck. Moreover, suppose that $G_\kappa \geq G_\Diamond$ at each point on this neck. If $\varepsilon \leq \varepsilon_\Diamond$, then $G_\kappa$ is pointwise non-decreasing under surgery. Moreover, if $\varepsilon \leq \varepsilon_\Diamond$ and $\sigma \leq \sigma_\Diamond$, then, for each $\delta \geq 0$, the quantity 
\[\max \Big \{ G_\kappa^{\sigma-1} \, \Big ( H - \frac{(n-1)^2(n+2)}{4} \, (1+\delta) \, G_\kappa \Big ),0 \Big \}\] 
is pointwise non-increasing under surgery.
\end{lemma}

Note that the constants $G_\Diamond,\varepsilon_\Diamond,\sigma_\Diamond$ in Lemma \ref{properties.of.surgery} do not depend on $\delta$. \\

\textbf{Proof.} 
We argue as in the proof of Theorem 5.3 (ii) in \cite{Huisken-Sinestrari2} (see pp.~179--180 in that paper). As in \cite{Huisken-Sinestrari2}, we put $\Lambda=10$. Moreover, we define $u(z) = r_0 \, e^{-\frac{B}{z-\Lambda}}$, where $r_0$ denotes the radius of the neck on which we perform surgery and $B$ is a large positive constant which will be specified later.

We first consider the region $S^{n-1} \times (\Lambda,3\Lambda]$. Let $\lambda_1,\hdots,\lambda_n$ denote the curvature eigenvalues at a point on the original neck, and let $\tilde{\lambda}_1,\hdots,\tilde{\lambda}_n$ be the curvature eigenvalues at the corresponding point on the bent hypersurface. Given any $\theta>0$, we can choose the parameter $B$ and the curvature scale $G_\Diamond$ sufficiently large so that 
\[|\tilde{\lambda}_1 - (\lambda_1+\tau_0 \, D_1 D_1 u+\tau_0 u \lambda_1^2)| \leq \theta \tau_0 \, D_1 D_1 u\] 
and 
\[|\tilde{\lambda}_i - (\lambda_i+\tau_0 u \lambda_i^2)| \leq \theta \tau_0 \, D_1 D_1 u\] 
for $i=2,\hdots,n$ (cf. \cite{Huisken-Sinestrari2}, (3.38)). Moreover, by choosing $B$ sufficiently large, we can arrange that $u \leq \theta \, r_0^2 \, D_1 D_1 u$ and $D_1 D_1 u \leq \theta \, r_0^{-1}$ (see \cite{Huisken-Sinestrari2}, Lemma 3.18). Therefore, we obtain 
\[|\tilde{\lambda}_1 - (\lambda_1+\tau_0 \, D_1 D_1 u)| \leq \theta \tau_0 \, D_1 D_1 u\] 
and 
\[|\tilde{\lambda}_i - \lambda_i| \leq \theta \tau_0 \, D_1 D_1 u\] 
for $i=2,\hdots,n$. This implies 
\[\tilde{H} = H + \tau_0 \, D_1 D_1 u + O(\theta \tau_0 \, D_1 D_1 u)\] 
and 
\[\tilde{G}_\kappa^{-1} = G_\kappa^{-1} - \tau_0 \, D_1 D_1 u \, \sum_{1<j} \frac{1}{(\lambda_1+\lambda_j-2\kappa)^2} + O(\theta \, r_0^2 \, \tau_0 \, D_1 D_1 u).\] 
In particular, if we choose $B$ sufficiently large, then we have $\tilde{G}_\kappa \geq G_\kappa$. We next compute 
\[\tilde{G}_\kappa^{-1} \tilde{H} = G_\kappa^{-1} H - \tau_0 \, D_1 D_1 u \, \Big ( \sum_{1<j} \frac{H}{(\lambda_1+\lambda_j-2\kappa)^2} - G_\kappa^{-1} \Big ) + O(\theta \, r_0 \, \tau_0 \, D_1 D_1 u).\] 
On an exact cylinder, we have $\sum_{1<j} \frac{H}{(\lambda_1+\lambda_j)^2} = \frac{4(n-1)}{n+2} \, G^{-1}$. Hence, if $\varepsilon$ and $r_0$ are sufficiently small, then we have $\sum_{1<j} \frac{H}{(\lambda_1+\lambda_j-2\kappa)^2} \geq \frac{3n-2}{n+2} \, G_\kappa^{-1}$ at each point on the original neck. Hence, if we choose $B$ large enough, then we obtain 
\[\tilde{G}_\kappa^{-1} \tilde{H} \leq G_\kappa^{-1} H - \frac{n-2}{n+2} \, G_\kappa^{-1} \, \tau_0 \, D_1 D_1 u.\] From this, we deduce  
\[\tilde{G}_\kappa^{\sigma-1} \tilde{H} \leq G_\kappa^{\sigma-1} H - \frac{n-2}{n+2} \, G_\kappa^{\sigma-1} \, \tau_0 \, D_1 D_1 u + O(\sigma \, G_\kappa^{\sigma-1} \, \tau_0 \, D_1 D_1 u).\] 
Hence, if $\sigma$ is sufficiently small, then we have 
\[\tilde{G}_\kappa^{\sigma-1} \tilde{H} \leq G_\kappa^{\sigma-1} H.\] 
Since $\tilde{G}_\kappa \geq G_\kappa$, we obtain 
\begin{align*} 
&\max \Big \{ \tilde{G}_\kappa^{\sigma-1} \, \Big ( \tilde{H} - \frac{(n-1)^2(n+2)}{4} \, (1+\delta) \, \tilde{G}_\kappa \Big ),0 \Big \} \\ 
&\leq \max \Big \{ G_\kappa^{\sigma-1} \, \Big ( H - \frac{(n-1)^2(n+2)}{4} \, (1+\delta) \, G_\kappa \Big ),0 \Big \} 
\end{align*}
at each point in the region $S^{n-1} \times (\Lambda,3\Lambda]$ and for each $\delta \geq 0$. 

Finally, we consider the region $S^{n-1} \times [3\Lambda,4\Lambda]$. Having fixed the surgery parameters $B,\tau_0$, we can choose $\varepsilon$ sufficiently small so that 
\[\tilde{H} - \frac{(n-1)^2(n+2)}{4} \, \tilde{G}_\kappa \leq 0\] 
in the region $S^{n-1} \times [3\Lambda,4\Lambda]$. Therefore, for each $\delta \geq 0$, we have 
\[\max \Big \{ \tilde{G}_\kappa^{\sigma-1} \, \Big ( \tilde{H} - \frac{(n-1)^2(n+2)}{4} \, (1+\delta) \, \tilde{G}_\kappa \Big ),0 \Big \} = 0\] 
in the region $S^{n-1} \times [3\Lambda,4\Lambda]$. This completes the proof of Lemma \ref{properties.of.surgery}. \\

From now on, we will assume that the surgery parameters are chosen so that Lemma \ref{properties.of.surgery} applies. \\

\begin{lemma}
\label{lower.bound.for.G.surgery.version}
Suppose that $M_t$, $t \in [0,T)$, is a surgically modified flow starting from a closed, embedded, $\kappa$-two-convex hypersurface $M_0$. If the curvature tensor of the ambient manifold satisfies $\overline{R}_{1313}+\overline{R}_{2323} \geq -2\kappa^2$ at each point on $M_t$, then $\inf_{M_t} G_\kappa$ approaches infinity in finite time.
\end{lemma} 

\textbf{Proof.} 
In between surgery times, we have 
\[\frac{\partial}{\partial t} G_\kappa \geq \frac{\partial G_\kappa}{\partial h_{ij}} \, D_i D_j G_\kappa + \frac{1}{C(n)} \, G_\kappa^3.\] 
We claim that $\inf_{M_t} G_\kappa$ is non-decreasing across each surgery time. To see this, suppose that $t$ is a surgery time, and that $x \in M_{t+}$ is a point in the surgically modified region. By Lemma \ref{properties.of.surgery}, there exists a point $y \in M_{t-}$ such that $G_\kappa(x,t+) \geq G_\kappa(y,t-)$. Consequently, $\inf_{M_{t+}} G_\kappa \geq \inf_{M_{t-}} G_\kappa$. From this, the assertion follows easily. \\

\begin{proposition} 
\label{uniformly.two.convex.surgery.version}
Suppose that $M_t$, $t \in [0,T)$, is a surgically modified flow starting from a closed, embedded, $\kappa$-two-convex hypersurface $M_0$. Then there exists a uniform constant $\beta$, depending only on $T$, $\kappa$, the initial hypersurface $M_0$, and the ambient manifold, such that $G_\kappa \geq \beta H$ at each point on $M_t$.
\end{proposition}

\textbf{Proof.} 
In between surgery times, we have 
\[\frac{\partial}{\partial t} H \leq \frac{\partial G_\kappa}{\partial h_{ij}} \, (D_i D_j H + h_{ik} \, h_{jk} \, H) + C \, H.\] 
Moreover, the ratio $\frac{H}{G_\kappa}$ is uniformly bounded from above in the surgery regions. Hence, the maximum principle implies that the ratio $\frac{H}{G_\kappa}$ is uniformly bounded from above on bounded time intervals. This completes the proof of Proposition \ref{uniformly.two.convex.surgery.version}. \\

As above, Proposition \ref{uniformly.two.convex.surgery.version} implies that $\frac{\partial G_\kappa}{\partial h_{ij}} \geq \beta^2 \, g_{ij}$. In particular, the equation is uniformly parabolic. 

We next establish a cylindrical estimate for surgically modified flows.

\begin{proposition} 
\label{cylindrical.estimate.surgery.version}
Suppose that $M_t$, $t \in [0,T)$, is a surgically modified flow starting from a closed, embedded, $\kappa$-two-convex hypersurface $M_0$. Let $\delta$ be an arbitrary positive real number. Then 
\[H \leq \frac{(n-1)^2(n+2)}{4} \, (1+\delta) \, G_\kappa + C,\] 
where $C$ is a positive constant that depends only on $\delta$, $T$, $\kappa$, the initial hypersurface $M_0$, and the ambient manifold.
\end{proposition}

\textbf{Proof.} 
Let $\sigma_\Diamond$ be as in Lemma \ref{properties.of.surgery}. For $\sigma \leq \sigma_\Diamond$, we define 
\[f_\sigma = G_\kappa^{\sigma-1} \, \Big ( H - \frac{(n-1)^2(n+2)}{4} \, (1+\delta) \, G_\kappa \Big ).\] 
In between surgery times, we have 
\[\frac{d}{dt} \bigg ( \int_{M_t} f_{\sigma,+}^p \bigg ) \leq (Cp)^p \, |M_t|,\] 
provided that $p$ is sufficiently large and $p^{\frac{1}{2}} \, \sigma$ is sufficiently small. Moreover, Lemma \ref{properties.of.surgery} guarantees that $\int_{M_t} f_{\sigma,+}^p$ is non-increasing across each surgery time. Hence, we can find a small positive constant $c_0$ such that 
\[\int_{M_t} f_{\sigma,+}^p \leq C\] 
for $p \geq \frac{1}{c_0}$ and $\sigma \leq c_0 \, p^{-\frac{1}{2}}$, where $C$ is a positive constant that depends only on $p$, $\sigma$, $\delta$, $T$, $\kappa$, the initial hypersurface $M_0$, and the ambient manifold. 

In the next step, we fix $p$ and $\sigma$ such that $p$ is very large and $0 < \sigma < c_0 \, (2np)^{-\frac{1}{2}} - 2 \, p^{-1}$. Moreover, let 
\[f_{\sigma,k} = G_\kappa^{\sigma-1} \, \Big ( H - \frac{(n-1)^2(n+2)}{4} \, (1+\delta) \, G_\kappa \Big ) - k\] 
and 
\[f_{\sigma,k,+} = \max \{f_{\sigma,k},0\}.\] 
Finally, we put $A(k) = \int_0^T |M_t \cap \{f_{\sigma,k} \geq 0\}|$. Then 
\begin{align*}
\frac{d}{dt} \bigg ( \int_{M_t} f_{\sigma,k,+}^p \bigg )
&\leq -\frac{1}{C} \int_{M_t} f_{\sigma,k,+}^{p-2} \, |\nabla f_{\sigma,k}|^2 \\
&+ C \int_{M_t} f_{\sigma,k,+}^{p-1} \, f_\sigma \, G_\kappa^2 + C \, |M_t \cap \{f_{\sigma,k} \geq 0\}|. 
\end{align*} 
Moreover, by Lemma \ref{properties.of.surgery}, the quantity $\int_{M_t} f_{\sigma,k,+}^p$ is non-increasing across each surgery time. This implies 
\[\sup_{t \in [0,T)} \int_{M_t} f_{\sigma,k,+}^p \leq C \, A(k) + C \int_0^T \int_{M_t} G_\kappa^2 \, f_{\sigma,k,+}^{p-1} \, f_\sigma\] 
and 
\[\int_0^T \int_{M_t} f_{\sigma,k,+}^{p-2} \, |\nabla f_{\sigma,k}|^2 \leq C \, A(k) + C \int_0^T \int_{M_t} G_\kappa^2 \, f_{\sigma,k,+}^{p-1} \, f_\sigma\]
for $k$ sufficiently large. Here, $C$ is a positive constant independent of $k$. Arguing as in the smooth case, we obtain 
\[A(\tilde{k})^{1-\frac{1}{n+1}} \, (\tilde{k}-k)^p \leq C \, A(k)^{1-\frac{1}{2n}}\] 
provided that $\tilde{k} \geq k$ and $k$ is sufficiently large. Here, $C$ is a positive constant independent of $k$ and $\tilde{k}$. Iterating this inequality gives $A(k)=0$ for some positive constant $k = k(p,\sigma,\delta,T,\kappa,M_0,X)$. Thus, $f_\sigma \leq k$ everywhere. This completes the proof of Proposition \ref{cylindrical.estimate.surgery.version}. \\

Combining Proposition \ref{cylindrical.estimate.surgery.version} with Proposition \ref{alg}, we can draw the following conclusion:

\begin{corollary}
\label{convexity.estimate.surgery.version}
Suppose that $M_t$, $t \in [0,T)$, is a surgically modified flow starting from a closed, embedded, $\kappa$-two-convex hypersurface $M_0$, and let $\delta$ be an arbitrary positive real number. Then 
\[\lambda_1 \geq -\delta \, G_\kappa - C,\] 
where $C$ is a positive constant that depends only on $\delta$, $T$, $\kappa$, the initial hypersurface $M_0$, and the ambient manifold.
\end{corollary}

In the next step, we verify that the inscribed radius estimate remains valid for surgically modified flows.

\begin{proposition} 
\label{noncollapsing.surgery.version}
Suppose that $M_t$, $t \in [0,T)$, is a surgically modified flow starting from a closed, embedded, $\kappa$-two-convex hypersurface $M_0$. Then the inscribed radius is bounded from below by $\frac{\alpha}{G_\kappa}$ at each point on $M_t$. Here, $\alpha$ is a positive constant that depends only on $T$, $\kappa$, the initial hypersurface $M_0$, and the ambient manifold.
\end{proposition}

\textbf{Proof.} 
Let $\mu$ be the quantity introduced in Section \ref{inscr.rad}. In between surgery times, we have 
\[\frac{\partial}{\partial t} \mu \leq \sum_{i,j} \frac{\partial G_\kappa}{\partial h_{ij}} \, (D_i D_j \mu + h_{ik} \, h_{jk} \, \mu) + C \, \mu\] 
whenever $\frac{\mu}{G_\kappa}$ is sufficiently large. 

In the next step, we claim that the ratio $\frac{\mu}{G_\kappa}$ is uniformly bounded from above in the surgery regions. To see this, suppose that $t$ is a surgery time and $N \subset M_{t-}$ is a neck on which we perform surgery.  Then the interior of the solid tube associated with $N$ is disjoint from $M_{t-}$ (see \cite{Huisken-Sinestrari2}, Theorem 3.26). Consequently, the ratio $\frac{\mu}{G_\kappa}$ is uniformly bounded from above on the neck $N$, and also on the cap which is inserted during surgery. 

Using the maximum principle, we conclude that the ratio $\frac{\mu}{G_\kappa}$ is uniformly bounded from above on bounded time intervals. Since the inscribed radius is comparable to $\frac{1}{\mu}$, the assertion follows. \\

Our next goal is to establish a pointwise curvature derivative estimate for surgically modified flows. We begin by extending the curvature estimates for radial graphs to the case of flows with surgery.

\begin{lemma}
\label{estimate.for.surgery.scale}
There exists a positive real number $\Xi \geq 100$, depending only on $n$, with the following property. Let $r \leq 1$ and let $\Omega_t$, $t \in [-r^2,0]$, be a one-parameter family of smooth open domains such that $B_r(p) \subset \Omega_t$ and the hypersurfaces $M_t = \partial \Omega_t$ form a surgically modified flow with surgery scale $G_*$. Moreover, suppose that $\langle -\exp_x^{-1}(p),\nu \rangle \geq 10^{-3} \, r$ at each point $x \in \partial U_t \cap B_{2r}(p)$, where $U_t$ denotes the connected component of $\Omega_t \cap B_{2r}(p)$ which contains the ball $B_r(p)$. If $G_* r \geq \Xi$, then the set $\partial U_t \cap B_{\frac{5r}{3}}(p)$ is free of surgeries for each $t \in [-r^2,0]$.
\end{lemma}

\textbf{Proof.} 
Suppose that the set $\partial U_{t+} \cap B_{\frac{5r}{3}}(p)$ contains a point modified by surgery. If $G_* r$ is sufficiently large, then the hypersurface $\partial U_{t+} \cap B_{2r}(p)$ contains an $(\varepsilon,6,10)$-neck. Moreover, if $G_* r$ is sufficiently large, this neck violates the star-shapedness condition $\langle -\exp_x^{-1}(p),\nu \rangle \geq 10^{-3} \, r$. Thus, we conclude that $G_* r$ is bounded from above by a large constant, as claimed. \\

\begin{proposition} 
\label{interior.curvature.estimate.for.radial.graphs.surgery.version}
Let $r \leq 1$ and let $\Omega_t$, $t \in [-r^2,0]$, be a one-parameter family of smooth open domains such that $B_r(p) \subset \Omega_t$ and the hypersurfaces $M_t = \partial \Omega_t$ form a surgically modified flow with surgery scale $G_*$. Moreover, we assume that $\langle -\exp_x^{-1}(p),\nu \rangle \geq 10^{-3} \, r$ and $G_\kappa \geq \beta H$ at each point $x \in \partial U_t \cap B_{2r}(p)$, where $U_t$ denotes the connected component of $\Omega_t \cap B_{2r}(p)$ which contains the ball $B_r(p)$. Then the norm of the second fundamental form satisfies 
\[\Big ( \frac{25r^2}{9}-d(p,x)^2 \Big ) \, (t+r^2)^{\frac{1}{2}} \, |h| \leq C \, r^2\] 
for all $t \in [-r^2,0]$ and all $x \in \partial U_t \cap B_{\frac{5r}{3}}(p)$. Here, $C$ is a positive constant that depends only on $\beta$ and the ambient manifold.
\end{proposition}

\textbf{Proof.} 
Let 
\[\eta = \frac{25r^2}{9}-d(p,x)^2.\] 
Moreover, we define 
\[\psi = \Big ( \frac{25r^2}{9}-d(p,x)^2 \Big ) \, (\langle -\exp_x^{-1}(p),\nu \rangle^2 - 10^{-7} \, r^2)^{-\frac{1}{2}} \, G_\kappa\] 
and 
\[Q(t) = \sup_{x \in \partial U_t \cap B_{\frac{5r}{3}}(p)} \psi(x,t).\] 
Let us fix a time $t \in [-r^2,0]$, and let $x \in \partial U_t \cap B_{\frac{5r}{3}}(p)$ be a point satisfying $\psi(x,t) = Q(t)$. Suppose first that $x$ lies in the surgically modified region. In this case, $G_\kappa(x,t) \leq C \, G_*$. Moreover, Lemma \ref{estimate.for.surgery.scale} implies that $G_* r \leq \Xi$. Putting these facts together, we conclude that $G_\kappa(x,t) \leq C \, r^{-1}$. This gives $Q(t) = \psi(x,t) \leq C$. 

Consequently, if $\psi(x,t) = Q(t)$ and $Q(t)$ is sufficiently large, then $x$ does not lie in the surgically modified region. In particular, we have $Q(t-) \geq Q(t+)$ if $Q(t+)$ is sufficiently large. Arguing as in the proof of Proposition \ref{interior.curvature.estimate.for.radial.graphs}, we conclude that 
\[\frac{\partial}{\partial t} \psi \leq -\frac{1}{C} \, r^2 \, \eta^{-2} \, \psi^3 \leq -\frac{1}{C} \, r^{-2} \, \psi^3,\] 
provided that $\psi(x,t) = Q(t)$ and $Q(t)$ is sufficiently large. Hence, if $Q(t)$ is sufficiently large, then we have 
\[\limsup_{t' \nearrow t} \frac{Q(t)-Q(t')}{t-t'} \leq -\frac{1}{C} \, r^{-2} \, Q(t)^3.\] 
This finally gives 
\[Q(t) \leq C \, r \, (t+r^2)^{-\frac{1}{2}}\] 
for all $t \in [-r^2,0]$. Since $|h| \leq C \, G_\kappa$, the assertion follows. \\

\begin{corollary} 
\label{higher.regularity.for.radial.graphs.surgery.version}
Let $r \leq 1$ and let $\Omega_t$, $t \in [-r^2,0]$, be a one-parameter family of smooth open domains such that $B_r(p) \subset \Omega_t$ and the hypersurfaces $M_t = \partial \Omega_t$ form a surgically modified flow with surgery scale $G_*$. Finally, we assume that $\langle -\exp_x^{-1}(p),\nu \rangle \geq 10^{-3} \, r$ and $G_\kappa \geq \beta H$ at each point $x \in \partial U_t \cap B_{2r}(p)$, where $U_t$ denotes the connected component of $\Omega_t \cap B_{2r}(p)$ which contains the ball $B_r(p)$. Then 
\[r^2 \, |\nabla h(x,0)| + r^3 \, |\nabla^2 h(x,0)| \leq \Lambda\] 
for all points $x \in \partial U_0 \cap B_{\frac{4r}{3}}(p)$ satisfying $G_\kappa(x,0) \geq \alpha \, r^{-1}$. Here, $\Lambda$ is a positive constant that depends only on $\alpha$, $\beta$, and the ambient manifold.
\end{corollary}

\textbf{Proof.} 
By Proposition \ref{interior.curvature.estimate.for.radial.graphs.surgery.version}, we can find a positive constant $K \geq 100$ such that $|h| \leq K \, r^{-1}$ and $G_\kappa \leq K \, r^{-1}$ for all $t \in [-\frac{r^2}{4},0]$ and all $x \in \partial U_t \cap B_{\frac{3r}{2}}(p)$. Let us fix an arbitrary point $x \in \partial U_0 \cap B_{\frac{4r}{3}}(p)$ satisfying $G_\kappa(x,0) \geq \alpha \, r^{-1}$. Let $\Xi$ be the constant in Lemma \ref{estimate.for.surgery.scale}, and let $\tau \in [-\frac{r^2}{100K^2},0]$ be the smallest number with the property that the parabolic neighborhood $\mathcal{P}(x,0,\frac{r}{4\Xi},|\tau|)$ is free of surgeries. For each $t \in (\tau,0]$, we have $\mathcal{P}(x,0,\frac{r}{4\Xi},|\tau|) \cap M_t \subset \partial U_t \cap B_{\frac{3r}{2}}(p)$. In particular, we have $|h| \leq K \, r^{-1}$ and $G_\kappa \leq K \, r^{-1}$ at each point in $\mathcal{P}(x,0,\frac{r}{4\Xi},|\tau|)$. Hence, on the set $\mathcal{P}(x,0,\frac{r}{4\Xi},|\tau|)$, the function $G_\kappa$ satisfies a uniformly parabolic equation with bounded coefficients. We now distinguish two cases: 

\textit{Case 1:} Suppose first that $\tau=-\frac{r^2}{100K^2}$. Using the Krylov-Safonov theorem (cf. Theorem \ref{krylov.safonov.interior.estimate}), we obtain a H\"older estimate for the function $G_\kappa$ on the set $\mathcal{P}(x,0,\frac{r}{8\Xi},\frac{r^2}{200K^2})$. In particular, there exists a uniform constant $\theta \in (0,\min \{\frac{1}{16\Xi},\frac{1}{400K^2}\})$ such that $\frac{\alpha}{2} \, r^{-1} \leq G_\kappa \leq K \, r^{-1}$ on the set $\mathcal{P}(x,0,\theta \, r,\theta \, r^2)$. Theorem \ref{interior.estimate.for.fully.nonlinear.PDE} gives H\"older estimates for the second fundamental form on the set $\mathcal{P}(x,0,\frac{\theta}{2} \, r,\frac{\theta}{2} \, r^2)$. Using Schauder theory, we obtain estimates for all derivatives of the second fundamental form on the set $\mathcal{P}(x,0,\frac{\theta}{4} \, r,\frac{\theta}{4} \, r^2)$. In particular, this gives bounds for $|\nabla h(x,0)|$ and $|\nabla^2 h(x,0)|$.

\textit{Case 2:} Suppose next that $\tau > -\frac{r^2}{100K^2}$. In this case, the set  $\mathcal{P}(x,0,\frac{r}{4\Xi},|\tau|) \cap M_{\tau+}$ contains a point $q$ which lies in a surgery region. Since $q \in \partial U_{\tau+} \cap B_{\frac{3r}{2}}(p)$, we have $G_* r \leq \Xi$ by Lemma \ref{estimate.for.surgery.scale}. Moreover, we have 
\[G_*^{-2} \, |\nabla h| + G_*^{-3} \, |\nabla^2 h| + G_*^{-4} \, |\nabla^3 h| \leq C(n)\] 
for all points $p \in M_{\tau+}$ satisfying $d_{M_{\tau+}}(p,q) \leq G_*^{-1}$. Finally, using the inequalities $|h| \leq K \, r^{-1}$ and $G_\kappa \leq K \, r^{-1}$, we conclude that the intrinsic diameter of the set $\mathcal{P}(x,0,\frac{r}{4\Xi},|\tau|) \cap M_{\tau+}$ is bounded from above by $\frac{r}{\Xi}$. This implies that 
\begin{align*} 
\mathcal{P}(x,0,\frac{r}{4\Xi},|\tau|) \cap M_{\tau+} 
&\subset \{p \in M_{\tau+}: d_{M_{\tau+}}(p,q) \leq \frac{r}{\Xi}\} \\ 
&\subset \{p \in M_{\tau+}: d_{M_{\tau+}}(p,q) \leq G_*^{-1}\}. 
\end{align*}
Hence, we obtain 
\begin{align*} 
&\frac{r^2}{\Xi^2} \, |\nabla h| + \frac{r^3}{\Xi^3} \, |\nabla^2 h| + \frac{r^4}{\Xi^4} \, |\nabla^3 h| \\ 
&\leq G_*^{-2} \, |\nabla h| + G_*^{-3} \, |\nabla^2 h| + G_*^{-4} \, |\nabla^3 h| \leq C(n) 
\end{align*}
on the set $\mathcal{P}(x,0,\frac{r}{4\Xi},|\tau|) \cap M_{\tau+}$. A version of the Krylov-Safonov theorem (cf. Corollary \ref{krylov.safonov.estimate.up.to.time.0}) now gives a H\"older bound for $G_\kappa$ on the set $\mathcal{P}(x,0,\frac{r}{8\Xi},|\tau|)$. In particular, there exists a uniform constant $\theta \in (0,\frac{1}{16\Xi})$ such that $\frac{\alpha}{2} \, r^{-1} \leq G_\kappa \leq K \, r^{-1}$ on the set $\mathcal{P}(x,0,\theta \, r,\min \{\theta \, r^2,|\tau|\})$. Corollary \ref{estimate.for.fully.nonlinear.PDE.up.to.time.0} now gives H\"older estimates for the second fundamental form on the set $\mathcal{P}(x,0,\frac{\theta}{2} \, r,\min \{\frac{\theta}{2} \, r^2,|\tau|\})$. Using Schauder theory, we obtain estimates for the first and second derivatives of the second fundamental form on the set $\mathcal{P}(x,0,\frac{\theta}{4} \, r,\min \{\frac{\theta}{4} \, r^2,|\tau|\})$. In particular, this gives an upper bound for $|\nabla h(x,0)|$ and $|\nabla^2 h(x,0)|$. This completes the proof of Corollary \ref{higher.regularity.for.radial.graphs.surgery.version}. \\

We are now in a position to prove a pointwise curvature derivative estimate for surgically modified flows.

\begin{theorem} 
\label{curvature.derivative.estimate.surgery.version}
Let us fix a closed, embedded, $\kappa$-two-convex hypersurface $M_0 = \partial \Omega_0$ in a Riemannian manifold, and a real number $\kappa \geq 0$. We can find a constant $G_\#$, depending only on $\kappa$, $M_0$, and the ambient manifold, such that the following holds. Suppose that $\Omega_t$, $t \in [0,T)$, is a one-parameter family of smooth open domains with the property that the hypersurfaces $M_t = \partial \Omega_t$ form a surgically modified flow starting from $M_0$ with surgery scale $G_* \geq G_\#$. Then we have 
\[\alpha^2 \, G_\kappa^{-2} \, |\nabla h| + \alpha^3 \, G_\kappa^{-3} \, |\nabla^2 h| \leq \Lambda\] 
for all points in spacetime satisfying $G_\kappa \geq G_\#$. Here, $\alpha$ is the constant in Proposition \ref{noncollapsing.surgery.version}, and $\Lambda$ is the constant appearing in Corollary \ref{higher.regularity.for.radial.graphs.surgery.version}.
\end{theorem} 

\textbf{Proof.} 
Suppose that the assertion is false. Then there exists a sequence of surgically modified flows $\mathcal{M}^{(k)}$ with surgery scales $G_*^{(k)} \to \infty$, and a sequence of points $(x_k,t_k) \in \mathcal{M}^{(k)}$ such that $G_\kappa(x_k,t_k) \to \infty$ and 
\[\alpha^2 \, G_\kappa(x_k,t_k)^{-2} \, |\nabla h(x_k,t_k)| + \alpha^3 \, G_\kappa(x_k,t_k)^{-3} \, |\nabla^2 h(x_k,t_k)| > \Lambda\] 
for each $k$. Using a standard point-picking argument, we can find, for each $k$, a point $(\bar{x}_k,\bar{t}_k) \in \mathcal{M}^{(k)}$ with the following properties: 
\begin{itemize}
\item[(i)] $\bar{t}_k \leq t_k$. 
\item[(ii)] $G_\kappa(\bar{x}_k,\bar{t}_k) \geq G_\kappa(x_k,t_k)$.
\item[(iii)] $\alpha^2 \, G_\kappa(\bar{x}_k,\bar{t}_k)^{-2} \, |\nabla h(\bar{x}_k,\bar{t}_k)| + \alpha^3 \, G_\kappa(\bar{x}_k,\bar{t}_k)^{-3} \, |\nabla^2 h(\bar{x}_k,\bar{t}_k)| > \Lambda$.
\item[(iv)] $\alpha^2 \, G_\kappa(x,t)^{-2} \, |\nabla h(x,t)| + \alpha^3 \, G_\kappa(x,t)^{-3} \, |\nabla^2 h(x,t)| \leq \Lambda$ for all points $(x,t) \in \mathcal{M}^{(k)}$ with $t \leq \bar{t}_k$ and $G_\kappa(x,t) \geq 2 \, G_\kappa(\bar{x}_k,\bar{t}_k)$. 
\end{itemize} 
For abbreviation, let $r_k = \alpha \, G_\kappa(\bar{x}_k,\bar{t}_k)^{-1}$. Note that $r_k \to 0$ in view of property (ii). Using Proposition \ref{noncollapsing.surgery.version}, we can find a point $p_k$ such that $\bar{x}_k \in \partial B_{r_k}(p_k)$ and $B_{r_k}(p_k) \subset \Omega_{\bar{t}_k}^{(k)}$. Clearly, $B_{r_k}(p_k) \subset \Omega_t^{(k)}$ for all $t \in [\bar{t}_k-r_k^2,\bar{t}_k]$. For each $t \in [\bar{t}_k-r_k^2,\bar{t}_k]$, we denote by $U_t^{(k)}$ the connected component of $\Omega_t^{(k)} \cap B_{2r_k}(p_k)$ which contains the ball $B_{r_k}(p_k)$. Clearly, the sets $U_t^{(k)}$ shrink as $t$ increases. We distinguish two cases:

\textit{Case 1:} Suppose that $C_{p_k,x} \subset U_t^{(k)}$ for all $t \in [\bar{t}_k-r_k^2,\bar{t}_k]$ and all points $x \in U_t^{(k)}$. This implies $\varangle (-\exp_x^{-1}(p_k),\nu) \leq \frac{\pi}{2}-\frac{1}{100}$ for all $t \in [\bar{t}_k-r_k^2,\bar{t}_k]$ and all points $x \in \partial U_t^{(k)} \cap B_{2r_k}(p_k)$. Consequently, $\langle -\exp_x^{-1}(p_k),\nu \rangle \geq 10^{-3} \, r_k$ for all $t \in [\bar{t}_k-r_k^2,\bar{t}_k]$ and all points $x \in \partial U_t^{(k)} \cap B_{2r_k}(p_k)$. Corollary \ref{higher.regularity.for.radial.graphs.surgery.version} gives $r_k^2 \, |\nabla h(x,\bar{t}_k)| + r_k^3 \, |\nabla^2 h(x,\bar{t}_k)| \leq \Lambda$ for all points $x \in \partial U_{\bar{t}_k}^{(k)} \cap B_{\frac{4r_k}{3}}(p_k)$ satisfying $G_\kappa(x,\bar{t}_k) \geq \alpha \, r_k^{-1}$. On the other hand, we clearly have $\bar{x}_k \in \partial U_{\bar{t}_k}^{(k)} \cap B_{\frac{4r_k}{3}}(p_k)$, $G_\kappa(\bar{x}_k,\bar{t}_k) = \alpha \, r_k^{-1}$, and furthermore $r_k^2 \, |\nabla h(\bar{x}_k,\bar{t}_k)| + r_k^3 \, |\nabla^2 h(\bar{x}_k,\bar{t}_k)| > \Lambda$ in view of property (iii) above. This is a contradiction.

\textit{Case 2:} Suppose that there exists a time $\tilde{t}_k \in [\bar{t}_k-r_k^2,\bar{t}_k]$ and a point $x \in U_{\tilde{t}_k}^{(k)}$ such that $C_{p_k,x} \not\subset U_{\tilde{t}_k}^{(k)}$. Let 
\[A^{(k)} = \{x \in U_{\tilde{t}_k}^{(k)}: C_{p_k,x} \subset U_{\tilde{t}_k}^{(k)}\}.\] 
It is clear that $B_{r_k}(p_k) \subset A^{(k)}$, $A^{(k)} \neq U_{\tilde{t}_k}^{(k)}$, and $A^{(k)}$ is relatively closed as a subset of $U_{\tilde{t}_k}^{(k)}$. Since $U_{\tilde{t}_k}^{(k)}$ is connected, it follows that $A^{(k)}$ cannot be an open set. Consequently, there exists a point $\tilde{x}_k \in A^{(k)}$ with the property that $B_\sigma(\tilde{x}_k) \not\subset A^{(k)}$ for all $\sigma>0$. Note that $r_k \leq |\tilde{x}_k-p_k| < 2r_k$ and $C_{p_k,\tilde{x}_k} \subset U_{\tilde{t}_k}^{(k)}$. Moreover, since 
\[\partial C_{p_k,\tilde{x}_k} \subset S_{p_k,\tilde{x}_k} \cup \{\tilde{x}_k\} \cup B_{\frac{r_k}{2}}(p_k) \subset S_{p_k,\tilde{x}_k} \cup U_{\tilde{t}_k}^{(k)},\] 
the hypersurface $S_{p_k,\tilde{x}_k}$ touches $\partial U_{\tilde{t}_k}^{(k)} \cap B_{2r_k}(p_k)$ somewhere from the inside. 

Let us consider a point $y_k$ where the hypersurface $S_{p_k,\tilde{x}_k}$ touches $\partial U_{\tilde{t}_k}^{(k)} \cap B_{2r_k}(p_k)$ from the inside. Clearly, $r_k \leq d(p_k,y_k) < 2r_k$. Since $C_{p_k,\tilde{x}_k}$ has aperture $2 \cdot \frac{1}{100}$, we can find a unit vector $v_k \in T_{y_k} X$ such that $\varangle (v_k,\nu(y_k,\tilde{t}_k)) \geq \frac{\pi}{2}+10^{-3}$ and 
\[\{\exp_{y_k}(v): 0 < |v| < \frac{r_k}{2}, \, \varangle (v_k,v) \leq 10^{-3}\} \subset C_{p_k,\tilde{x}_k} \subset U_{\tilde{t}_k}^{(k)}.\] 
In particular, we have 
\begin{equation} 
\label{c}
\{\exp_{y_k}(v): 0 < |v| < \frac{r_k}{2}, \, \varangle (v_k,v) \leq 10^{-3}\} \cap M_t^{(k)} = \emptyset. 
\end{equation}
for all $t \leq \tilde{t}_k$.

By Lemma \ref{smallest.curvature.eigenvalue.of.S}, the smallest curvature eigenvalue of $S_{p_k,\tilde{x}_k}$ is less than $-10^{-3} \, d(p_k,\tilde{x}_k)^{-1}$ at each point on $S_{p_k,\tilde{x}_k}$. Since the hypersurface $S_{p_k,\tilde{x}_k}$ touches $\partial U_{\tilde{t}_k}^{(k)} \cap B_{2r_k}(p_k)$ from the inside at $y_k$, it follows that 
\[\lambda_1(y_k,\tilde{t}_k) \leq -10^{-3} \, d(p_k,\tilde{x}_k)^{-1} \leq -\frac{10^{-3}}{2} \, r_k^{-1} = -\frac{10^{-3}}{2} \, \alpha^{-1} \, G_\kappa(\bar{x}_k,\bar{t}_k).\] 
In particular, $\lambda_1(y_k,\tilde{t}_k) \to -\infty$ in view of property (ii) above. Using Corollary \ref{convexity.estimate}, we obtain 
\[\lambda_1(y_k,\tilde{t}_k) \, G_\kappa(y_k,\tilde{t}_k)^{-1} \to 0.\] 
Thus, we conclude that 
\[G_\kappa(y_k,\tilde{t}_k) \, G_\kappa(\bar{x}_k,\bar{t}_k)^{-1} \to \infty.\] 
In particular, we have $G_\kappa(y_k,\tilde{t}_k) \geq 8 \, G_\kappa(\bar{x}_k,\bar{t}_k)$ if $k$ is sufficiently large. For each $k$, we define 
\[L_k = \min \Big \{ \inf \Big \{ G_\kappa(y_k,\tilde{t}_k) \, d_{M_{\tilde{t}_k}}(y_k,x): x \in M_{\tilde{t}_k}^{(k)}, \, \frac{G_\kappa(x,\tilde{t}_k)}{G_\kappa(y_k,\tilde{t}_k)} \notin [\frac{1}{2},2] \Big \},10^6 \Big \}.\] 
By definition of $L_k$, we have $\frac{1}{2} \, G_\kappa(y_k,\tilde{t}_k) \leq G_\kappa(x,\tilde{t}_k) \leq 2 \, G_\kappa(y_k,\tilde{t}_k)$ for all points $x \in M_{\tilde{t}_k}^{(k)}$ satisfying $d_{M_{\tilde{t}_k}}(y_k,x) \leq L_k \, G_\kappa(y_k,\tilde{t}_k)^{-1}$. Using property (iv) above, we can find a uniform constant $\theta \in (0,1)$ such that  
\[\sup_{\mathcal{P}_{\mathcal{M}^{(k)}}(y_k,\tilde{t}_k,(L_k+\theta) \, G_\kappa(y_k,\tilde{t}_k)^{-1},\tilde{\theta} \, G_\kappa(y_k,\tilde{t}_k)^{-2})} G_\kappa \leq 4 \, G_\kappa(y_k,\tilde{t}_k)\] 
and 
\[\inf_{\mathcal{P}_{\mathcal{M}^{(k)}}(y_k,\tilde{t}_k,(L_k+\theta) \, G_\kappa(y_k,\tilde{t}_k)^{-1},\tilde{\theta} \, G_\kappa(y_k,\tilde{t}_k)^{-2})} G_\kappa \geq \frac{1}{4} \, G_\kappa(y_k,\tilde{t}_k),\] 
whenever $\tilde{\theta} \in [0,\theta]$ and $\mathcal{P}_{\mathcal{M}^{(k)}}(y_k,\tilde{t}_k,(L_k+\theta) \, G_\kappa(y_k,\tilde{t}_k)^{-1},\tilde{\theta} \, G_\kappa(y_k,\tilde{t}_k)^{-2})$ is free of surgeries. 

For each $k$, we denote by $\tilde{\theta}_k \in [0,\theta]$ the largest number with the property that the parabolic neighborhood $\mathcal{P}_{\mathcal{M}^{(k)}}(y_k,\tilde{t}_k,(L_k+\theta) \, G_\kappa(y_k,\tilde{t}_k)^{-1},\tilde{\theta}_k \, G_\kappa(y_k,\tilde{t}_k)^{-2})$ is free of surgeries. We distinguish two subcases: 

\textit{Subcase 2.1.} Suppose that $\lim_{k \to \infty} \tilde{\theta}_k > 0$. In this case, we will argue that $y_k$ lies at the center of an $(\varepsilon_k,6,\frac{(n-1)(n+2)}{4} \cdot 10^5)$-neck in $M_{\tilde{t}_k}^{(k)}$ for some sequence $\varepsilon_k \to 0$. To prove this, we restrict the flow $\mathcal{M}^{(k)}$ to the parabolic neighborhood $\mathcal{P}_{\mathcal{M}^{(k)}}(y_k,\tilde{t}_k,(L_k+\theta) \, G_\kappa(y_k,\tilde{t}_k)^{-1},\tilde{\theta}_k \, G_\kappa(y_k,\tilde{t}_k)^{-2})$. On this parabolic neighborhood, the ratio $\frac{H}{G_\kappa(y_k,\tilde{t}_k)}$ is uniformly bounded from above, and the ratio $\frac{\lambda_1+\lambda_2-2\kappa}{G_\kappa(y_k,\tilde{t}_k)}$ is uniformly bounded from below. Hence, if we perform a parabolic dilation around the point $(y_k,\tilde{t}_k)$ with factor $G_\kappa(y_k,\tilde{t}_k)$, then the rescaled flow has bounded curvature and is uniformly two-convex. By property (iv) above, we have bounds for the first and second derivatives of the second fundamental form. Hence, the rescaled flows converge to a smooth, non-flat limit flow in $\mathbb{R}^{n+1}$, which moves with velocity $G$ and satisfies the pointwise inequality $H \leq \frac{(n-1)^2(n+2)}{4} \, G$ (see Theorem \ref{cylindrical.estimate.surgery.version}). Since $\lambda_1(y_k,\tilde{t}_k) < 0$ for each $k$, there exists a point on the limit flow where the smallest curvature eigenvalue is non-positive. Using Proposition \ref{splitting}, we conclude that the limit flow is contained in a family of shrinking cylinders. In particular, this implies 
\[\sup \{G_\kappa(x,\tilde{t}_k): x \in M_{\tilde{t}_k}^{(k)}, \, d_{M_{\tilde{t}_k}}(y_k,x) \leq L_k \, G_\kappa(y_k,\tilde{t}_k)^{-1}\} \leq (1+o(1)) \, G_\kappa(y_k,\tilde{t}_k)\] 
and 
\[\inf \{G_\kappa(x,\tilde{t}_k): x \in M_{\tilde{t}_k}^{(k)}, \, d_{M_{\tilde{t}_k}}(y_k,x) \leq L_k \, G_\kappa(y_k,\tilde{t}_k)^{-1}\} \geq (1-o(1)) \, G_\kappa(y_k,\tilde{t}_k).\] 
Consequently, we have $L_k=10^6$ if $k$ is sufficiently large. Moreover, the point $y_k$ lies at the center of an $(\varepsilon_k,6,\frac{(n-1)(n+2)}{4} \cdot 10^5)$-neck in $M_{\tilde{t}_k}^{(k)}$ for some sequence $\varepsilon_k \to 0$. Since $\varangle (v_k,\nu(y_k,\tilde{t}_k)) \geq \frac{\pi}{2}+10^{-3}$, we conclude that 
\begin{equation} 
\label{d}
\{\exp_{y_k}(sv_k): 0 < s < 10^4 \, G_\kappa(y_k,\tilde{t}_k)^{-1}\} \cap M_{\tilde{t}_k}^{(k)} \neq \emptyset 
\end{equation} 
if $k$ is sufficiently large. Since $G_\kappa(y_k,\tilde{t}_k) \, r_k \to \infty$, the statements (\ref{c}) and (\ref{d}) are in contradiction. 

\textit{Subcase 2.2.} Suppose finally that $\lim_{k \to \infty} \tilde{\theta}_k = 0$. Let $\hat{t}_k = \tilde{t}_k - \tilde{\theta}_k \, G_\kappa(y_k,\tilde{t}_k)^{-2}$. By definition of $\tilde{\theta}_k$, the set $\mathcal{P}_{\mathcal{M}^{(k)}}(y_k,\tilde{t}_k,(L_k+\theta) \, G_\kappa(y_k,\tilde{t}_k)^{-1},\tilde{\theta}_k \, G_\kappa(y_k,\tilde{t}_k)^{-2}) \cap M_{\hat{t}_k+}^{(k)}$ contains a point $q_k$ which lies in a surgery region. Then the hypersurface $\{x \in M_{\hat{t}_k+}^{(k)}: d_{M_{\hat{t}_k+}}(q_k,x) \leq 10^8 \, (G_*^{(k)})^{-1}\}$ is a capped-off neck. Moreover, we have $\frac{1}{2} \, G_*^{(k)} \leq G_\kappa(q_k,\hat{t}_k+) \leq 100 \, G_*^{(k)}$. On the other hand, we have $\frac{1}{4} \, G_\kappa(y_k,\tilde{t}_k) \leq G_\kappa(q_k,\hat{t}_k+) \leq 4 \, G_\kappa(y_k,\tilde{t}_k)$. Putting these facts together gives $\frac{1}{8} \, G_*^{(k)} \leq G_\kappa(y_k,\tilde{t}_k) \leq 400 \, G_*^{(k)}$. 

By following the point $y_k \in M_{\tilde{t}_k}^{(k)}$ backwards in time, we obtain a point $z_k \in M_{\hat{t}_k+}^{(k)}$ satisfying $d(y_k,z_k) \leq o(1) \, G_\kappa(y_k,\tilde{t}_k)^{-1} \leq o(1) \, (G_*^{(k)})^{-1}$ and $|\nu(y_k,\tilde{t}_k)-P_{z_k}^{y_k} \nu(z_k,\hat{t}_k+)| \leq o(1)$. Here, $P_{z_k}^{y_k}: T_{z_k} X \to T_{y_k} X$ denotes the parallel transport along a minimizing geodesic in $X$ starting at $z_k$ and ending at $y_k$. Note that $q_k$ and $z_k$ belong to the set $\mathcal{P}_{\mathcal{M}^{(k)}}(y_k,\tilde{t}_k,(L_k+\theta) \, G_\kappa(y_k,\tilde{t}_k)^{-1},\tilde{\theta}_k \, G_\kappa(y_k,\tilde{t}_k)^{-2}) \cap M_{\hat{t}_k+}^{(k)}$, and the intrinsic diameter of that set is at most $2 \, (L_k+\theta+o(1)) \, G_\kappa(y_k,\tilde{t}_k)^{-1}$. This gives 
\begin{align*} 
d_{M_{\hat{t}_k+}}(q_k,z_k) 
&\leq 2 \, (L_k+\theta+o(1)) \, G_\kappa(y_k,\tilde{t}_k)^{-1} \\ 
&\leq 2 \, (10^6+1) \, G_\kappa(y_k,\tilde{t}_k)^{-1} \\ 
&\leq 16 \, (10^6+1) \, (G_*^{(k)})^{-1} 
\end{align*} 
for $k$ large. In particular, if $k$ is sufficiently large, the point $z_k$ lies on the capped-off neck described above. Since $\varangle (v_k,\nu(y_k,\tilde{t}_k)) \geq \frac{\pi}{2}+10^{-3}$, we have $\varangle (v_k,P_{z_k}^{y_k} \nu(z_k,\hat{t}_k+)) \geq \frac{\pi}{2}+10^{-3}-o(1)$. This implies 
\begin{equation} 
\label{e}
\{\exp_{y_k}(v): 0 < |v| < 10^4 \, (G_*^{(k)})^{-1}, \, \varangle (v_k,v) \leq 10^{-3}\} \cap M_{\hat{t}_k+}^{(k)} \neq \emptyset 
\end{equation} 
if $k$ is sufficiently large. Since $G_*^{(k)} \, r_k \geq \frac{1}{400} \, G(y_k,\tilde{t}_k) \, r_k \to \infty$, the statement (\ref{e}) contradicts (\ref{c}). This completes the proof of Theorem \ref{curvature.derivative.estimate.surgery.version}. \\

The following result is the analogue of the Neck Detection Lemma in \cite{Huisken-Sinestrari2}.

\begin{theorem}[Neck Detection Lemma]
\label{neck.detection}
Let us fix closed embedded, $\kappa$-two-convex hypersurface $M_0$ in a Riemannian manifold. Given positive real numbers $\varepsilon_0,L_0,\theta>0$, we can find positive numbers $\eta_0,G_0>0$ with the following property: Let $M_t$, $t \in [0,T)$, be a surgically modified flow in $\mathbb{R}^{n+1}$ starting from $M_0$ with surgery scale $G_*$. Moreover, suppose that $t_0 \in [0,T)$ and $p_0 \in M_{t_0}$ satisfy 
\begin{itemize}
\item $G_\kappa(p_0,t_0) \geq G_0$, $\frac{\lambda_1(p_0,t_0)}{G_\kappa(p_0,t_0)} \leq \eta_0$, 
\item the parabolic neighborhood $\mathcal{P}(p_0,t_0,(L_0+1) \, G_\kappa(p_0,t_0)^{-1},\theta \, G_\kappa(p_0,t_0)^{-2})$ does not contain surgeries.
\end{itemize}
Then $(p_0,t_0)$ lies at the center of an $(\varepsilon_0,6,\frac{(n-1)(n+2)}{4} \, L_0)$-neck. \end{theorem}

Note that the constants $\eta_0$ and $G_0$ may depend on $\varepsilon_0,L_0,\theta$, $\kappa$, the initial hypersurface $M_0$, and the ambient manifold, but they are independent of the surgery parameters $\varepsilon,L$. \\

\textbf{Proof.} 
Suppose that the assertion is false. Then there exists a sequence of surgically modified flows $\mathcal{M}^{(k)}$ and a sequence of points $(p_k,t_k)$ with the following properties: 
\begin{itemize} 
\item $G_\kappa(p_k,t_k) \geq k$ and $\frac{\lambda_1(p_k,t_k)}{G_\kappa(p_k,t_k)} \leq \frac{1}{k}$. 
\item The parabolic neighborhood $\mathcal{P}_{\mathcal{M}^{(k)}}(p_k,t_k,(L_0+1) \, G_\kappa(p_k,t_k)^{-1},\theta \, G_\kappa(p_k,t_k)^{-2})$ is free of surgeries. 
\item The point $p_k$ does not lie at the center of an $(\varepsilon_0,6,\frac{(n-1)(n+2)}{4} \, L_0)$-neck in the hypersurface $M_{t_k}^{(k)}$.
\end{itemize}
For each $k$, we put 
\[L_k = \min \Big \{ \inf \Big \{ G_\kappa(p_k,t_k) \, d_{M_{t_k}}(p_k,x): x \in M_{t_k}^{(k)}, \, \frac{G_\kappa(x,t_k)}{G_\kappa(p_k,t_k)} \notin [\frac{1}{2},2] \Big \},L_0 \Big \}.\] 
By definition of $L_k$, we have $\frac{1}{2} \, G_\kappa(p_k,t_k) \leq G_\kappa(x,t_k) \leq 2 \, G_\kappa(p_k,t_k)$ for all points $x \in M_{t_k}^{(k)}$ satisfying $d_{M_{t_k}}(p_k,x) \leq L_k \, G_\kappa(p_k,t_k)^{-1}$. Using Theorem \ref{curvature.derivative.estimate.surgery.version}, we obtain 
\[\sup_{\mathcal{P}_{\mathcal{M}^{(k)}}(p_k,t_k,(L_k+\tilde{\theta}) \, G_\kappa(p_k,t_k)^{-1},\tilde{\theta} \, G_\kappa(p_k,t_k)^{-2})} G_\kappa \leq 4 \, G_\kappa(p_k,t_k)\] 
and 
\[\inf_{\mathcal{P}_{\mathcal{M}^{(k)}}(p_k,t_k,(L_k+\tilde{\theta}) \, G_\kappa(p_k,t_k)^{-1},\tilde{\theta} \, G_\kappa(p_k,t_k)^{-2})} G_\kappa \geq \frac{1}{4} \, G_\kappa(p_k,t_k)\] 
for some uniform constant $\tilde{\theta} \in (0,\theta)$ which is independent of $k$. 

In the next step, we restrict the flow $\mathcal{M}^{(k)}$ to the parabolic neighborhood $\mathcal{P}_{\mathcal{M}^{(k)}}(p_k,t_k,(L_k+\tilde{\theta}) \, G_\kappa(p_k,t_k)^{-1},\tilde{\theta} \, G_\kappa(p_k,t_k)^{-2})$. On this parabolic neighborhood, the ratio $\frac{H}{G_\kappa(p_k,t_k)}$ is uniformly bounded from above, and the ratio $\frac{\lambda_1+\lambda_2-2\kappa}{G_\kappa(p_k,t_k)}$ is uniformly bounded from below. Hence, if we perform a parabolic dilation around the point $(p_k,t_k)$ with factor $G_\kappa(p_k,t_k)$, then the rescaled flow has bounded curvature and is uniformly two-convex. By Theorem \ref{curvature.derivative.estimate.surgery.version}, the first and second derivatives of the second fundamental form are bounded as well. Hence, the rescaled flows converge to a smooth, non-flat limit flow in $\mathbb{R}^{n+1}$, which moves with normal velocity $G$ and satisfies the pointwise inequality $H \leq \frac{(n-1)^2(n+2)}{4} \, G$ (see Theorem \ref{cylindrical.estimate.surgery.version}). Since $\frac{\lambda_1(p_k,t_k)}{G_\kappa(p_k,t_k)} \leq \frac{1}{k}$ for each $k$, there exists a point on the limit flow where the smallest curvature eigenvalue is non-positive. Again, Proposition \ref{splitting} implies that the limit flow is contained in a family of shrinking cylinders. In particular, 
\[\sup \{G_\kappa(x,t_k): x \in M_{t_k}^{(k)}, \, d_{M_{t_k}}(p_k,x) \leq L_k \, G_\kappa(p_k,t_k)^{-1}\} \leq (1+o(1)) \, G_\kappa(p_k,t_k)\] 
and 
\[\inf \{G_\kappa(x,t_k): x \in M_{t_k}^{(k)}, \, d_{M_{t_k}}(p_k,x) \leq L_k \, G_\kappa(p_k,t_k)^{-1}\} \geq (1-o(1)) \, G_\kappa(p_k,t_k).\] 
Thus, we conclude that $L_k=L_0$ if $k$ is sufficiently large. Moreover, the point $p_k$ lies at the center of an $(\varepsilon_k,6,\frac{(n-1)(n+2)}{4} \cdot L_0)$-neck in $M_{t_k}^{(k)}$ for some sequence $\varepsilon_k \to 0$. This is a contradiction. \\

\section{Existence of surgically modified flows}

\label{surgery.riemannian.case}

In this final section, we outline how we can implement the surgery algorithm from \cite{Huisken-Sinestrari2}. We first consider the case that the ambient manifold is the Euclidean space $\mathbb{R}^{n+1}$. Having established the convexity estimate, the cylindrical estimate, and the curvature derivative estimate for surgically modified flows, the arguments in Section 7 and Section 8 of \cite{Huisken-Sinestrari2} carry over unchanged to our setting. Thus, we can use the surgery algorithm in \cite{Huisken-Sinestrari2} to extend the flow beyond singularities. This proves the assertion in the special case when the ambient space is the Euclidean space $\mathbb{R}^{n+1}$.

In the remainder of this section, we sketch how the arguments in Section 7 and Section 8 of \cite{Huisken-Sinestrari2} can be adapted to the Riemannian setting. Let us fix an ambient Riemannian manifold $X$. We assume that the surgery parameters are chosen as explained on pp.~208--209 of \cite{Huisken-Sinestrari2}. This fixes the values of all surgery parameters except the curvature threshold $H_1$, which we may choose arbitrarily large (cf. the remark at the bottom of p.~209 in \cite{Huisken-Sinestrari2}). The basic idea is that, by choosing $H_1$ sufficiently large, the curvature of the background metric becomes negligible and will not interfere with the proof of the Neck Continuation Theorem. There are only two points in the proof of the Neck Continuation Theorem that require minor modifications: 

First, in the proof of the Neck Continuation Theorem on p.~214, one considers a unit vector field $\omega$ in ambient space. One then considers the flow on $M_{t_0}$ generated by the vector field $\frac{\omega^T}{|\omega^T|^2}$, where $\omega^T$ denotes the projection of $\omega$ to the tangent space of $M_{t_0}$ (see \cite{Huisken-Sinestrari2}, p.~205). In the Euclidean setting, $\omega$ is parallel, and consequently we have $\frac{d}{dy} \langle \nu,\omega \rangle \geq \lambda_1$ along each trajectory of this ODE, where $\lambda_1$ denotes the smallest eigenvalue of the second fundamental form (cf. \cite{Huisken-Sinestrari2}, Proposition 7.18). In the Riemannian setting, we choose a local height function $u$ in ambient space such that $|\nabla u|=1$ at each point where $u$ is defined. Note that $u$ is defined on a small geodesic ball in ambient space; the radius of that ball is a positive constant which depends only on the ambient manifold $X$. We then consider the flow on $M_{t_0}$ generated by the vector field $\frac{\omega^T}{|\omega^T|^2}$, where $\omega = \nabla u$ and $\omega^T$ denotes the projection of $\omega$ to the tangent space of $M_{t_0}$. Along each trajectory of the ODE, we have 
\begin{align*} 
\frac{d}{dy} \langle \nu,\omega \rangle 
&= \frac{\langle \bar{D}_{\omega^T} \nu,\omega \rangle + \langle \nu,\bar{D}_{\omega^T} \omega \rangle}{|\omega^T|^2} \\ 
&= \frac{h(\omega^T,\omega^T) + \langle \nu - \langle \nu,\omega \rangle \, \omega,\bar{D}_{\omega^T} \omega \rangle}{|\omega^T|^2}, 
\end{align*}
where in the last step we have used that $\omega$ has unit length. Using the identity $|\nu-\langle \nu,\omega \rangle \, \omega| = |\omega-\langle \omega,\nu \rangle \, \nu| = |\omega^T|$, we conclude that 
\[\frac{d}{dy} \langle \nu,\omega \rangle \geq \lambda_1 - |\bar{D} \omega|.\] 
The error term $|\bar{D} \omega|$ does not affect the proof of the Neck Continuation Theorem on pp.~214--215 of \cite{Huisken-Sinestrari2}. Indeed, in the region $z \in [\bar{z},z^*]$, the smallest eigenvalue of the second fundamental form is bounded from below by $\lambda_1 \geq \eta_1 H$. Moreover, the mean curvature $H$ is larger than $\frac{H_1}{4}$ in the region $z \in [\bar{z},z^*]$ (cf. \cite{Huisken-Sinestrari2}, p.~211). Hence, the piece of the neck where $z \in [\bar{z},z^*]$ has diameter $O(H_1^{-1})$; in particular, if $H_1$ is sufficiently large, then the piece of the neck where $z \in [\bar{z},z^*]$ is contained in the domain of definition of the height function $u$. Furthermore, we have $\frac{d}{dy} \langle \nu,\omega \rangle \geq \lambda_1 - |\bar{D} \omega| \geq \frac{\eta_1 H_1}{4} - |\bar{D} \omega|$ for $z \in [\bar{z},z^*]$. Since $\eta_1$ has already been chosen at this stage, we can now choose the curvature threshold $H_1$ sufficiently large so that $\frac{d}{dy} \langle \nu,\omega \rangle > 0$ for $z \in [\bar{z},z^*]$, which is all we need for the argument to work. 

Second, on p.~215, one needs to show that the part of the surface coming after $\Sigma_{y'}$ is a convex cap. To that end, one again considers the flow on $M_{t_0}$ generated by the vector field $\frac{\omega^T}{|\omega^T|^2}$. In the Euclidean setting, one can show that the inequalities 
\[\langle \nu,\omega \rangle < 1, \qquad \lambda_1 > 0, \qquad H > \frac{H_1}{4\Theta}, \qquad \langle \nu,\omega \rangle > \varepsilon_1\] 
hold for all $y \in [y',y_{\text{\rm max}})$. This argument requires a minor modification in the Riemannian setting. To explain this, let $\eta_2$ be the constant introduced in the third application of the Neck Detection Lemma (see \cite{Huisken-Sinestrari2}, p.~209, statement (P5)). We claim that the inequalities 
\begin{equation} 
\tag{$\star$}
\langle \nu,\omega \rangle < 1, \qquad \lambda_1 \geq \eta_2 H, \qquad H > \frac{H_1}{4\Theta}, \qquad \langle \nu,\omega \rangle > \varepsilon_1
\end{equation}
hold for all $y \in [y',y_{\text{\rm max}})$, provided that the curvature threshold $H_1$ is chosen sufficiently large. Indeed, the inequalities in ($\star$) are clearly satisfied for $y=y'$. If one of the inequalities in ($\star$) fails for some $y>0$, we consider the smallest value of $y$ for which that happens. The first inequality in ($\star$) cannot fail first by definition of $y_{\text{\rm max}}$. If the second inequality in ($\star$) is the first one to fail, then we have $\lambda_1 = \eta_2 H$ at some point on that slice. In view of our choice of $\eta_2$, we conclude that this point lies on a cylindrical graph of length $5$ and $C^1$-norm less than $\varepsilon_1$ (see \cite{Huisken-Sinestrari2}, p.~209, statement (P5)), but this is ruled out by the fourth inequality in ($\star$). If the third inequality in ($\star$) is the first one to fail, we obtain a contradiction with Lemma 7.19 in \cite{Huisken-Sinestrari2}. Finally, as long as ($\star$) holds, we have $\frac{d}{dy} \langle \nu,\omega \rangle \geq \lambda_1 - |\bar{D} \omega| \geq \eta_2 H - |\bar{D} \omega| \geq \frac{\eta_2 H_1}{4\Theta} - |\bar{D} \omega|$. Note that $\eta_2$ and $\Theta$ have already been fixed at this stage. Hence, if we choose the curvature threshold $H_1$ sufficiently large, then $\langle \nu,\omega \rangle$ is montone increasing along each trajectory of the ODE. This implies that the fourth inequality in ($\star$) cannot fail first. Thus, the inequalities in ($\star$) hold for all $y \in [y',y_{\text{\rm max}})$. Consequently, the union of the surfaces $\Sigma_y$ is a convex cap, which is precisely what we need to complete the proof of the Neck Continuation Theorem. This completes our discussion of the Riemannian case.

Finally, if the curvature tensor of the ambient manifold satisfies $\overline{R}_{1313}+\overline{R}_{2323} \geq -2\kappa^2$ at each point in $\Omega_0$, then Lemma \ref{lower.bound.for.G.surgery.version} guarantees that the flow becomes extinct in finite time. This completes the proof of Theorem \ref{existence.of.surgically.modified.flows}.

\appendix 

\section{Review of Krylov-Safonov estimates}

For the convenience of the reader, we collect some well-known regularity results for parabolic equations. The first one is the crucial H\"older estimate of Krylov and Safonov \cite{Krylov-Safonov} (see also \cite{Krylov2}, Theorem 7 on pp.~137--138):

\begin{theorem}[N.V.~Krylov, M.V.~Safonov]
\label{krylov.safonov.interior.estimate}
Let $v: B_1(0) \times [0,1] \to \mathbb{R}$ be a solution of the parabolic equation $\frac{\partial}{\partial t} v = \sum_{i,j} a_{ij} \, D_i D_j v + \sum_i b_i \, D_i v + f$. We assume that the coefficients satisfy $\frac{1}{K} \, \delta_{ij} \leq a_{ij} \leq K \, \delta_{ij}$ and $|b_i| \leq K$. Then  
\[[v]_{C^{\gamma;\frac{\gamma}{2}}(B_{\frac{1}{2}}(0) \times [\frac{1}{2},1])} \leq C \, \Big ( \sup_{B_1(0) \times [0,1])} v - \inf_{B_1(0) \times [0,1]} v \Big ) + C \, \|f\|_{C^0(B_1(0) \times [0,1])},\] 
where $\gamma>0$ and $C>0$ depend only on $K$.
\end{theorem}

\begin{corollary}[N.V.~Krylov, M.V.~Safonov]
\label{krylov.safonov.estimate.up.to.time.0}
Let $0 < \tau \leq \frac{1}{4}$, and let $v: B_1(0) \times [0,\tau] \to \mathbb{R}$ be a solution of the parabolic equation $\frac{\partial}{\partial t} v = \sum_{i,j} a_{ij} \, D_i D_j v + \sum_i b_i \, D_i v + f$. We assume that the coefficients satisfy $\frac{1}{K} \, \delta_{ij} \leq a_{ij} \leq K \, \delta_{ij}$ and $|b_i| \leq K$. Finally, we assume that $\|v\|_{C^0(B_1(0) \times [0,\tau])} + \|v(\cdot,0)\|_{C^2(B_1(0))} + \|f\|_{C^0(B_1(0) \times [0,\tau])} \leq L$. Then  
\[[v]_{C^{\gamma;\frac{\gamma}{2}}(B_{\frac{1}{2}}(0) \times [0,\tau])} \leq C.\] 
Here, $\gamma>0$ depends only on $K$, and $C$ depends only on $K$ and $L$. In particular, $\gamma$ and $C$ are independent of $\tau$.
\end{corollary}

\textbf{Proof.} 
We sketch the argument for the convenience of the reader. Using a straightforward barrier argument, we can show that 
\[\sup_{B_r(x) \times [0,\min\{r^2,\tau\}]} v \leq v(x,0)+Cr\] 
and 
\[\inf_{B_r(x) \times [0,\min\{r^2,\tau\}]} v \geq v(x,0)-Cr\] 
for $x \in B_{\frac{1}{2}}(0)$ and $0 < r \leq \frac{1}{2}$. This gives 
\begin{equation} 
\label{oscillation.bound}
\sup_{B_r(x) \times [0,\min\{r^2,\tau\}]} v - \inf_{B_r(x) \times [0,\min\{r^2,\tau\}]} v \leq Cr 
\end{equation}
for $x \in B_{\frac{1}{2}}(0)$ and $0 < r \leq \frac{1}{2}$. Using Theorem \ref{krylov.safonov.interior.estimate} together with (\ref{oscillation.bound}), we obtain 
\begin{equation}
\label{holder.bound}
[v]_{C^{\gamma;\frac{\gamma}{2}}(B_{\frac{r}{2}}(x) \times [\frac{r^2}{2},r^2])} \leq C \, r^{1-\gamma} 
\end{equation}
for $x \in B_{\frac{1}{2}}(0)$ and $0 < r \leq \tau^{\frac{1}{2}}$. We now consider two points $(x,t)$ and $(\tilde{x},\tilde{t})$ in spacetime such that $t \geq \tilde{t} \geq 0$. If $2 \, |x-\tilde{x}| + 2 \, (t-\tilde{t})^{\frac{1}{2}} < t^{\frac{1}{2}}$, then (\ref{holder.bound}) gives 
\[|v(x,t) - v(\tilde{x},\tilde{t})| \leq C \, (|x-\tilde{x}| + (t-\tilde{t})^{\frac{1}{2}})^\gamma.\]  
On the other hand, if $2 \, |x-\tilde{x}| + 2 \, (t-\tilde{t})^{\frac{1}{2}} \geq t^{\frac{1}{2}}$, then (\ref{oscillation.bound}) implies 
\[|v(x,t) - v(\tilde{x},\tilde{t})| \leq C \, (|x-\tilde{x}| + t^{\frac{1}{2}}) \leq C \, (3 \, |x-\tilde{x}| + 2 \, (t-\tilde{t})^{\frac{1}{2}}).\] 
Putting these facts together, the assertion follows. \\

Combining the Krylov-Safonov estimate with the deep work of Evans \cite{Evans}, Krylov \cite{Krylov1}, and Caffarelli \cite{Caffarelli} on fully nonlinear elliptic equations gives:

\begin{theorem}
\label{interior.estimate.for.fully.nonlinear.PDE}
Let $u: B_1(0) \times [0,1] \to \mathbb{R}$ be a solution of a fully nonlinear parabolic equation 
\[\frac{\partial}{\partial t} u = \Phi(D^2 u,Du,u,x),\] 
where $\Phi$ depends smoothly on all its arguments. We assume that $u$ is bounded in $C^{2;1}(B_1(0) \times [0,1])$. Moreover, we assume that the equation is uniformly parabolic, and $\Phi$ is concave in the first argument. Then $u$ is uniformly bounded in $C^{2,\gamma;1,\frac{\gamma}{2}}(B_{\frac{1}{4}}(0) \times [\frac{1}{2},1])$ for some uniform constant $\gamma>0$.
\end{theorem}

\textbf{Proof.} 
Consider the function $v = \frac{\partial}{\partial t} u$. The function $v$ satisfies a uniformly parabolic equation. Moreover, $v$ is bounded in $C^0(B_1(0) \times [0,1])$, so the Krylov-Safonov estimate implies that $v$ is bounded in $C^{\gamma;\frac{\gamma}{2}}(B_{\frac{1}{2}}(0) \times [\frac{1}{2},1])$. Using Theorem 3 in \cite{Caffarelli}, it follows that $\sup_{t \in [\frac{1}{2},1]} \|u(t)\|_{C^{2,\gamma}(B_{\frac{1}{4}}(0))}$ is bounded from above. In other words, $D^2 u$ is uniformly H\"older continuous in space. Finally, 
\begin{align*} 
&\|D^2 u(t)-D^2 u(t')\|_{C^0(B_{\frac{1}{4}}(0))} \\ 
&\leq C \, \|D^2 u(t)-D^2 u(t')\|_{C^\gamma(B_{\frac{1}{4}}(0))}^{\frac{2}{2+\gamma}} \, \|u(t)-u(t')\|_{C^0(B_{\frac{1}{4}}(0))}^{\frac{\gamma}{2+\gamma}} \leq C \, |t-t'|^{\frac{\gamma}{2+\gamma}} 
\end{align*}
for $t,t' \in [\frac{1}{2},1]$. This shows that $D^2 u$ is uniformly H\"older continuous in time. \\

\begin{corollary}
\label{estimate.for.fully.nonlinear.PDE.up.to.time.0}
Let $0 < \tau \leq \frac{1}{4}$, and let $u: B_1(0) \times [0,\tau] \to \mathbb{R}$ be a solution of a fully nonlinear parabolic equation 
\[\frac{\partial}{\partial t} u = \Phi(D^2 u,Du,u,x),\] 
where $\Phi$ depends smoothly on all its arguments. We assume that $u$ is bounded in $C^{2;1}(B_1(0) \times [0,\tau])$, and that the initial function $u(\cdot,0)$ is bounded in $C^4(B_1(0))$. Moreover, we assume that the equation is uniformly parabolic, and $\Phi$ is concave in the first argument. Then $u$ is uniformly bounded in $C^{2,\gamma;1,\frac{\gamma}{2}}(B_{\frac{1}{4}}(0) \times [0,1])$ for some uniform constant $\gamma>0$.
\end{corollary}

\textbf{Proof.} 
We again consider the function $v = \frac{\partial}{\partial t} u$. The function $v$ satisfies a uniformly parabolic equation. Moreover, $v$ is bounded in $C^0(B_1(0) \times [0,1])$ and the initial function $v(\cdot,0)$ is bounded in $C^2(B_1(0))$. Consequently, $v$ is bounded in $C^{\gamma;\frac{\gamma}{2}}(B_{\frac{1}{2}}(0) \times [0,1])$. As above, Theorem 3 in \cite{Caffarelli} implies that $\sup_{t \in [0,1]} \|u(t)\|_{C^{2,\gamma}(B_{\frac{1}{4}}(0))}$ is uniformly bounded from above. In other words, $D^2 u$ is uniformly H\"older continuous in space. As above, we have the estimate 
\begin{align*} 
&\|D^2 u(t)-D^2 u(t')\|_{C^0(B_{\frac{1}{4}}(0))} \\ 
&\leq C \, \|D^2 u(t)-D^2 u(t')\|_{C^\gamma(B_{\frac{1}{4}}(0))}^{\frac{2}{2+\gamma}} \, \|u(t)-u(t')\|_{C^0(B_{\frac{1}{4}}(0))}^{\frac{\gamma}{2+\gamma}} \leq C \, |t-t'|^{\frac{\gamma}{2+\gamma}} 
\end{align*}
for $t,t' \in [0,1]$. Hence, $D^2 u$ is uniformly H\"older continuous in time, as claimed. \\


\begin{thebibliography}{99}
\bibitem{Andrews1}
B.~Andrews, \textit{Contraction of convex hypersurfaces in Riemannian spaces,} J. Diff. Geom. 39, 407--431 (1994)

\bibitem{Andrews2} 
B.~Andrews, \textit{Non-collapsing in mean-convex mean curvature flow,} Geom. Topol. 16, 1413--1418 (2012)

\bibitem{Andrews-Langford-McCoy1}
B.~Andrews, M.~Langford, and J.~McCoy, \textit{Non-collapsing in fully non-linear curvature flows,} Ann. Inst. H. Poincar\'e Anal. Non Lin\'eaire 30, 23--32 (2013)

\bibitem{Andrews-Langford-McCoy2}
B.~Andrews, M.~Langford, and J.~McCoy, \textit{Convexity estimates for hypersurfaces moving by convex curvature functions,} Anal. PDE 7, 407--433 (2014)

\bibitem{Ball} 
J.M.~Ball, \textit{Differentiability properties of symmetric and isotropic functions,} Duke Math. J. 51, 699--728 (1984)

\bibitem{Brendle1}
S.~Brendle, \textit{Embedded minimal tori in $S^3$ and the Lawson conjecture,} Acta Math. 211, 177--190 (2013)

\bibitem{Brendle2}
S.~Brendle, \textit{Two-point functions and their applications in geometry,} Bull. Amer. Math. Soc. 51, 581--596 (2014)

\bibitem{Brendle3}
S.~Brendle, \textit{A sharp bound for the inscribed radius under mean curvature flow,} Invent. Math. 202, 217--237 (2015)

\bibitem{Brendle4}
S.~Brendle, \textit{An inscribed radius estimate for mean curvature flow in Riemannian manifolds,} Ann. Scuola Norm. Sup. Pisa 16, 1447--1472 (2016)

\bibitem{Brendle-Huisken1}
S.~Brendle and G.~Huisken, \textit{Mean curvature flow with surgery of mean convex surfaces in $\mathbb{R}^3$,} Invent. Math. 203, 615--654 (2016)

\bibitem{Brendle-Huisken2}
S.~Brendle and G.~Huisken, \textit{Mean curvature flow with surgery of mean convex surfaces in three-manifolds,} to appear in J. Eur. Math. Soc.

\bibitem{Brendle-Hung}
S.~Brendle and P.K.~Hung, \textit{A sharp inscribed radius estimate for fully nonlinear flows,} to appear in Amer. J. Math.

\bibitem{Caffarelli}
L.A.~Caffarelli, \textit{Interior a-priori estimates for solutions of fully nonlinear equations,} Ann. of Math. 130, 189--213 (1989)

\bibitem{Ecker-Huisken1}
K.~Ecker and G.~Huisken, \textit{Mean Curvature Evolution of Entire Graphs,} Ann. of Math. 130, 453--471 (1989)

\bibitem{Ecker-Huisken2}
K.~Ecker and G.~Huisken, \textit{Interior estimates for hypersurfaces moving by mean curvature,} Invent. Math. 105, 547--569 (1991)

\bibitem{Evans}
L.C.~Evans, \textit{Classical solutions of fully nonlinear, convex, second-order elliptic equations,} Comm. Pure Appl. Math. 35, 333--363 (1982) 

\bibitem{Gerhardt}
C.~Gerhardt, \textit{Flow of nonconvex hypersurfaces into spheres,} J. Diff. Geom. 32, 299--314 (1990) 

\bibitem{Grayson}
M.~Grayson, \textit{The heat equation shrinks embedded plane curves to round points,} J. Diff. Geom. 26, 285--314 (1987)

\bibitem{Hamilton1} 
R.~Hamilton, \textit{Four-manifolds with positive curvature operator,} J. Diff. Geom. 24, 153--179 (1986)

\bibitem{Hamilton2}
R.~Hamilton, \textit{Isoperimetric estimates for the curve shrinking flow in the plane,} Modern Methods in Complex Analysis (Princeton 1992), 201--222, Ann. of Math. Stud. 137, Princeton University Press, Princeton NJ (1995)

\bibitem{Hamilton3}
R.~Hamilton, \textit{The formation of singularities in the Ricci flow,} Surveys in Differential Geometry, vol. II, 7--136, International Press, Somerville MA (1995)

\bibitem{Hamilton4} 
R.~Hamilton, \textit{Four-manifolds with positive isotropic curvature,} Comm. Anal. Geom. 5, 1--92 (1997)

\bibitem{Haslhofer-Kleiner1}
R.~Haslhofer and B.~Kleiner, \textit{Mean curvature flow of mean convex hypersurfaces,} to appear in Comm. Pure Appl. Math.

\bibitem{Haslhofer-Kleiner2}
R.~Haslhofer and B.~Kleiner, \textit{Mean curvature flow with surgery,} to appear in Duke Math. J.

\bibitem{Huisken}
G.~Huisken, \textit{A distance comparison principle for evolving curves,} Asian J. Math. 2, 127--133 (1998)

\bibitem{Huisken-Sinestrari1}
G.~Huisken and C.~Sinestrari, \textit{Convexity estimates for mean curvature flow and singularities of mean convex surfaces,} Acta Math. 183, 45--70 (1999)

\bibitem{Huisken-Sinestrari2}
G.~Huisken and C.~Sinestrari, \textit{Mean curvature flow with surgeries of two-convex hypersurfaces,} Invent. Math. 175, 137--221 (2009)

\bibitem{Krylov1} 
N.V.~Krylov, \textit{Boundedly inhomogeneous elliptic and parabolic equations,} Izv. Akad. Nauk SSSR Ser. Mat. 46, 487--523, 670 (1982)

\bibitem{Krylov2}
N.V.~Krylov, \textit{Nonlinear elliptic and parabolic equations of the second order,} 
Translated from the Russian by P.L.~Buzytsky, Mathematics and its Applications (Soviet Series), vol. 7, Reidel, Dordrecht (1987) 

\bibitem{Krylov-Safonov}
N.V.~Krylov and M.V.~Safonov, \textit{A property of the solutions of parabolic equations with measurable coefficients,} Izv. Akad. Nauk SSSR Ser. Mat. 44, 161–175, 239 (1980) 

\bibitem{Michael-Simon}
J.H.~Michael and L.M.~Simon, \textit{Sobolev and mean value inequalities on generalized submanifolds of $\mathbb{R}^n$,} Comm. Pure Appl. Math. 26, 316--379 (1973)

\bibitem{Perelman1} 
G.~Perelman, \textit{The entropy formula for the Ricci flow and its geometric applications,} arxiv:0211159

\bibitem{Perelman2}
G.~Perelman, \textit{Ricci flow with surgery on three-manifolds,} arxiv:0303109

\bibitem{Perelman3}
G.~Perelman, \textit{Finite extinction time for solutions to the Ricci flow on certain three-manifolds,} arxiv:0307245

\bibitem{Sha}
J.P.~Sha, \textit{$p$-convex Riemannian manifolds,} Invent. Math. 83, 437--447 (1986)

\bibitem{Sheng-Wang}
W.~Sheng and X.J.~Wang, \textit{Singularity profile in the mean curvature flow,} Methods Appl. Anal. 16, 139--155 (2009) 

\bibitem{Tso}
K.~Tso, \textit{Deforming a hypersurface by its Gauss-Kronecker curvature,} Comm. Pure Appl. Math. 38, 867--882 (1985)

\bibitem{Urbas}
J.~Urbas, \textit{On the expansion of starshaped hypersurfaces by symmetric functions of their principal curvatures,} Math. Z. 205, 355--372 (1990) 

\bibitem{White1}
B.~White, \textit{The size of the singular set in mean curvature flow of mean convex sets,} J. Amer. Math. Soc. 13, 665--695 (2000)

\bibitem{White2}
B.~White, \textit{The nature of singularities in mean curvature flow of mean convex sets,} J. Amer. Math. Soc. 16, 123--138 (2003)
\end{thebibliography}
\end{document}